\documentclass[16pt]{article}
\usepackage{amsmath,amsfonts,amssymb}
\usepackage{verbatim}
\usepackage{latexsym}
\newtheorem{theorem}{Theorem}[section]

\newtheorem{proposition}[theorem]{Proposition}
\newtheorem{definition}[theorem]{Definition}

\newtheorem{remark}[theorem]{Remark}

\newcommand\RR{{\Bbb R}}
\newcommand\CC{{\Bbb C}}

\newcommand\NN{{\Bbb N}}
\newcommand\ZZ{{\Bbb Z}}

\newcommand\diam{{\rm diam}}
\newcommand\supp{{\rm supp}}
\newcommand\sgn{{\rm sgn}}
\newcommand{\edth}{\partial}
\newcommand{\bedth}{\overline{\partial}}
\newcommand{\e}{\bf E}
\newcommand{\m}{\bf M}

\newcommand{\M}{\bf M}
\input{tcilatex}

\begin{document}
\title{Spin Wavelets on the Sphere}
\author{Daryl Geller\\
\footnotesize{Department of Mathematics, University of Stony Brook}\\
\footnotesize\texttt{{daryl@math.sunysb.edu}}\\ \\
Domenico Marinucci\\
\footnotesize{Department of Mathematics, University of Rome Tor Vergata}\\
\footnotesize\texttt{{marinucc@axp.mat.uniroma2.it}}}
\maketitle

\begin{abstract}
In recent years, a rapidly growing literature has focussed on the
construction of wavelet systems to analyze functions defined on the
sphere. Our purpose in this paper is to generalize these
constructions to situations where sections of line bundles,
rather than ordinary scalar-valued functions, are considered.
In particular, we propose {\em needlet-type spin wavelets} as an extension of the
needlet approach recently introduced by \cite{narc1,narc2} and then
considered for more general manifolds by \cite{gmcw,gmfr,gm3}. We discuss
localization properties in the real and harmonic domains, and
investigate stochastic properties for the analysis of spin random
fields. Our results are strongly motivated by cosmological
applications, in particular in connection to the analysis of Cosmic
Microwave Background polarization data.
\end{abstract}

\section{Introduction}

In recent years a rapidly growing literature has focussed on the
construction of wavelet systems on the sphere, see for instance
\cite{jfaa1}, \cite{jfaa3}  \cite{antoine,Antodem}  and the
references therein. These attempts have been motivated by strong
interest from the applied sciences, for instance in the areas of
Geophysics, Medical Imaging and especially Cosmology/Astrophysics.

As far as the latter are concerned, special emphasis has been
devoted to wavelet techniques for the statistical study of the
Cosmic Microwave Background (CMB) radiation data. These data can be
viewed as providing observations on the Universe in the immediate
sequel to the Big Bang, and as such they have been the object of
immense theoretical and applied interest over the last decade
\cite{dode2004}. In particular, experiments such as the WMAP
satellite from NASA have provided high resolution observations on
the "temperature" (i.e.\ intensity) of CMB radiation from all
directions of the full-sky (\cite{hinshaw}). These observations have
allowed precise estimates of several parameters of the greatest
interest for Cosmology and Theoretical Physics. Spherical wavelets
have found very extensive applications here, see for instance
(\cite{gorski,Cruz1,vielva04,chmpv,mcewen1,mcewen2,wiaux2}) and many
others. The rationale for such a widespread interest can be
explained as follows: CMB models are best analyzed in the frequency
domain, where the behavior at different multipoles can be
investigated separately; on the other hand, partial sky coverage and
other missing observations make the evaluation of exact spherical
harmonic transforms troublesome. The combination of these two
features makes the time-frequency localization properties of
wavelets most valuable.

Besides providing measurements on the radiation intensity,
experiments such as WMAP have also provided some preliminary
observations of a much more elusive physical entity, the so-called
polarization of the background radiation. Polarization is a property
of electromagnetic radiation whose physical significance is
described for instance in (\cite{kks96,zs}), see below for more
mathematical discussion. So far, empirical analysis of polarization
has been somewhat limited, because this signal is currently
measured with great difficulty. The situation with respect to CMB
polarization data, however, will significantly improve over the next
years, for instance by means of the ESA satellite mission Planck
(expected to be launched in Spring 2009) which will take full-sky
measurements of the polarized CMB sky (\cite{laureijs}) with much
greater precision. Moreover, both ESA and NASA are planning high
sensitivity full-sky satellite borne experiments within the next
10-20 years.

Polarization measurements are of extreme interest to physicists for
several reasons. Indeed, not only do they allow improved precision for
estimates of physical parameters which are already the focus of CMB
temperature data, but they also open entirely new areas of research.
Just to mention a striking example, at large scales the polarization
signal is expected to be dominated by a component related to a
background of gravitational waves which originated in the Big Bang
dynamics (in the so-called inflationary scenario, see
(\cite{dode2004}). Detection of this signal would be an outstanding
empirical validation of many Theoretical Cosmology claims, directly
related to Big Bang dynamic models but deeply rooted in General
Relativity. Taking into account the huge amount of polarization data
which will be available in the next 1-2 decades, as well as the
important cosmological information contained in these data, it is
clear that suitable mathematical tools for data analysis are in high
demand. While a large amount of mathematical statistics techniques
have been developed for analyzing CMB temperature data, far fewer
mathematical tools are available to analyze polarization data.

From the mathematical point of view, as we shall detail below,
polarization can be viewed as a (random) section of a line bundle on
the sphere. Our purpose here is then to extend spherical wavelet
constructions to the case where sections of line bundles, rather
than ordinary functions, are considered. To the best of our
knowledge, this is the first attempt to introduce wavelet techniques
for the case where one deals with more general mathematical structures
than ordinary (scalar-valued) functions on manifolds. We believe
this area can be expanded into disciplines other than cosmology - for
instance, tensor-valued random fields emerge naturally in brain
imaging data (\cite{armin}).

In particular, in this paper we shall extend to the line bundle case
the so-called needlet approach to spherical wavelets. Spherical
needlets were recently introduced by \cite{narc1,narc2}, and further
developed, and extended to general smooth compact Riemannian manifolds
in \cite{gmcw,gmfr,gm3}. In a random fields environment,
needlets were investigated by (\cite{BKMP06-3,bkmpBer}), with a view
to applications to the statistical analysis of CMB data;
applications in the physical literature include \cite{pbm06},
\cite{mpbb08}, \cite{guilloux}, see also (\cite{lan},
\cite{dela08,fay08,fg08}) and \cite{lanmar,mayun}.

More precisely, we will focus on \textit{needlet-type spin wavelets}, and we will
argue below that they enjoy properties analogous to those of the usual needlets
in the standard scalar case. In particular, we shall show below that
needlet-type spin wavelets enjoy both the localization and the uncorrelation
properties that make scalar needlets a powerful tool for the
analysis of scalar-valued spherical random fields. More details on
the plan and significance of our paper are given in the next
section.  

The companion article \cite{ghmkp}, intended primarily for physicists, goes into further
detail about the statistical uses of spin wavelets.  One should be aware, however, that there are some notational
differences between that article and this one.

\section{Plan of the Paper and Significance of the Results}

The spin $s$ concept that we will use in this paper originates in the fundamental work of
Newman and Penrose \cite{newpen}.
Writing for physicists, they say that a quantity $\eta$ defined on the sphere
has spin weight $s \in \ZZ$, provided that, whenever a tangent vector $m$
at a point transforms under coordinate change by $m' = e^{i\psi} m$, the quantity $\eta$, at that point, transforms by
$\eta' = e^{is\psi} \eta$.  They then develop a theory of spin $s$ spherical harmonics, upon which we
shall build to produce spin wavelets.  \\

In section 3, we put the notion of spin $s$ into acceptable mathematical language.
We do this in an elementary manner which remains very close to the spirit of the work of Newman and Penrose.
For another approach, using much more machinery, see \cite{eastod}.

Let ${\bf N}$ be the north pole $(0,0,1)$, let ${\bf S}$ be the south pole, $(0,0,-1)$, and let
$U_I$ be $S^2 \setminus \{{\bf N}, {\bf S}\}$.  If $R \in SO(3)$, we
define $U_R = RU_I$.  On $U_I$ we use standard spherical coordinates $(\theta,\phi)$
($0 < \theta < \pi$, $-\pi \leq \phi < \pi$), and analogously, on any $U_R$ we use coordinates
$(\theta_R,\phi_R)$ obtained by rotation of the coordinate system on $U_I$.
At each point $p$ of $U_R$ we let $\rho_R(p)$ be the unit tangent vector at $p$
which is tangent to the circle $\theta_R = $ constant, pointing in the direction of increasing $\phi_R$.
(Thus, for any $p \in U_I$, $\rho_R(Rp) = R_{*p} [\rho_{I}(p)]$.)   If $p \in U_{R_1} \cap U_{R_2}$,
we let $\psi_{p R_2 R_1}$ be the angle from $\rho_{R_1}(p)$ to $\rho_{R_2}(p)$.

Now say $\Omega \subseteq S^2$ is open.   We say that $f = (f_R)_{R \in SO(3)} \in C^{\infty}_s(\Omega)$, or
that $f$ is a smooth spin $s$ function on $\Omega$, provided all $f_R \in C^{\infty}(U_R \cap \Omega)$, and that
for all $R_1, R_2 \in SO(3)$ and all $p \in U_{R_1} \cap U_{R_2} \cap \Omega$,
\begin{equation}
\label{fr1r2}
f_{R_2}(p) = e^{is\psi}f_{R_1}(p),
\end{equation}
where $\psi = \psi_{p R_2 R_1}$.  If, say, $R_1 = I$, $R_2 = R$, then heuristically
$f_R$ is $f_I$ ``looked at after the coordinates have been rotated by $R$'';
at $p$, it has been multiplied by $e^{is\psi}$.   The angle $\psi$ would
be clearly be the same if we had instead chosen $\rho_R(p)$ to point in the direction of increasing $\theta_R$, say, so
this is an appropriate way to make Newman and Penrose's concept precise.  Note that if $s = 0$, we can clearly identify
$C^{\infty}_s(\Omega)$ with $C^{\infty}(\Omega)$.

For any $s$, we can identify $C^{\infty}_s(\Omega)$ with the
sections over $\Omega$ of the complex line bundle obtained by using the $e^{is\psi_{p R_2 R_1}}$ as transition functions from the
chart $U_{R_1}$ to the chart $U_{R_2}$; then $f_R$ is the trivialization of the section over $U_R \cap \Omega$.
This is the point of view that we take in section 3, since, as is often the case in mathematics, certain
properties are clearer, and more easily verified, if one uses the coordinate-free line bundle point of view.
However, right now, for the reader's convenience, we state our results without reference to line bundles.
Line bundles are rarely mentioned explicitly in this article after section 3.

We may define $L^2_s(\Omega)$
(resp.\ $C_s(\Omega)$), by requiring that the $f_R$ be in $L^2(U_R \cap \Omega)$
(resp.\ $C(U_R \cap \Omega)$) instead of $C^{\infty}(U_R \cap \Omega)$.
There is a well-defined inner product on $L^2_s(\Omega)$, given by $\langle f,g \rangle =
\langle f_R,g_R \rangle$; clearly this definition is independent of choice of $R$.
A key observation of section 3 is that there is a unitary action of $SO(3)$ on $L^2_s(S^2)$, given by $f \to f^R$, which is
determined by the equation $(f^R)_I(p) = f_R(Rp)$.  We think of $f^R$ as a ``rotate'' of $f$.\\

In section 4, we explain the spin $s$ theory of Newman and Penrose, within
the rigorous framework we have just outlined.  Most of our arguments are very close to theirs.

For $f$ smooth as above, following Newman and Penrose we define
\[ \edth_{sR} f_{R} = -(\sin \theta_R)^s \left(\frac{\partial}{\partial \theta_R} +
\frac{i}{\sin \theta_R}\frac{\partial}{\partial \phi_R} \right)(\sin \theta_R)^{-s} f_{R}, \]
and we show that the ``spin-raising'' operator
$\edth: C^{\infty}_s(\Omega) \to C^{\infty}_{s+1}(\Omega)$ given by
$(\edth f)_R = \edth_{sR}f_R$ is well-defined.  We also define $\bedth_{sR} f_R$ by
$\bedth_{sR} f_R = \overline{\edth_{sR} \overline{f_R}}$, which leads to the
spin-lowering operator $\bedth: C^{\infty}_s(\Omega) \to C^{\infty}_{s-1}(\Omega)$ given by
$(\bedth f)_R = \bedth_{sR}f_R$.   We show that $\edth$ and $\bedth$ commute with the
actions of $SO(3)$ on smooth spin functions.

Now, for $l \geq |s|$, let $b_{ls} = [(l+s)!/(l-s)!]^{1/2}$ for $s \geq 0$,
$b_{ls} = [(l-s)!/(l+s)!]^{1/2}$ for $s < 0$.  Let $\{Y_{lm}: l \geq 0,\: -l \leq m \leq l\}$ be the standard basis
of spherical harmonics
on $S^2$.  Following Newman and Penrose, for $l \geq |s|$, we define the spin $s$ spherical harmonics by
${}_sY_{lm} = \edth^s Y_{lm}/b_{ls}$ for $s \geq 0$, ${}_sY_{lm} = (-\bedth)^s Y_{lm}/b_{ls}$ for $s < 0$, and
we show that $\{{}_sY_{lm}: l \geq |s|, -l \leq m \leq l\}$ forms an orthonormal basis for
$L^2_s(S^2)$.  There is a (relatively) simple explicit expression for the $_sY_{lmI}$; it has the form
\begin{equation}
\label{sylmfrm}
_sY_{lmI}(\theta,\phi) =\ _{s}y_{lm}(\theta)e^{im\phi}
\end{equation}
for a suitable function $_sy_{lm}$.  (See (\ref{sphaspn}) for the explicit formula.)
 If $R \neq I$, $_sY_{lmR}$ does not have a simple
expression; but $_sY_{lmR}^R(\theta_R,\phi_R)$ equals $_sy_{lm}(\theta_R)e^{im\phi_R}$ in the coordinates on $U_R$.
(Here again $_sY_{lm}^R$ denotes the ``rotate'' of $_sY_{lm}$.)

There is also an analogue $\Delta_s$ of the spherical Laplacian for spin $s$ functions; specifically,
if $s \geq 0$, we let $\Delta_s = -\bedth\edth$, while if $s < 0$, we let $\Delta_s = -\edth\bedth$.
Then $\Delta_0$ is the usual spherical Laplacian.  For $l \geq |s|$, let
${\cal H}_{ls} = <{}_sY_{lm}: -l \leq m \leq l>$.   Also, for $l \geq |s|$, let
$\lambda_{ls} = (l-s)(l+s+1)$ if $s \geq 0$, and let $\lambda_{ls} = (l+s)(l-s+1)$ if $s < 0$.
Then ${\cal H}_{ls}$ is the subspace of $C^{\infty}_s(S^2)$ which consists of eigenfunctions of
$\Delta_s$ for the eigenvalue $\lambda_{ls}$.  This of course generalizes the situation in which $s = 0$.
All of the aforementioned facts follow by making the arguments of Newman and Penrose rigorous.\\

In sections 5, 6, 7 and 8, we present new results.
In section 5, we define operators with smooth kernels, from $C^{\infty}_s(S^2)$ to itself.
We say $K = (K_{R',R})_{R',R \in SO(3)}$
is such a smooth kernel if each\\ $K_{R',R} \in C^{\infty}(S^2 \times S^2)$, and if
$K_{{R_1}',{R_1}}(x,y) = e^{is(\psi'-\psi)}K_{R',R}(x,y)$
for all $x \in U_{{R_1}'} \cap U_{R'}$ and $y \in U_{{R_1}} \cap U_{R}$, where
$\psi' = \psi_{x R_1' R'}$ and $\psi = \psi_{y R R_1}$.  It is then evident that we may
consistently define an operator ${\cal K}: C^{\infty}_s(S^2) \to C^{\infty}_s(S^2)$ by
$({\cal K}f)_{R'}(x) = \int_{S^2} K_{R',R}(x,y) f_R(y) dS(y)$,
for all $x \in U_{R'}$; we call ${\cal K}$ the operator with kernel $K$.  (These operators may be
identified with operators with smooth kernels acting on sections of the aforementioned line bundle.)

Next we show that there is a good notion of spin $s$ {\em zonal harmonic}.  To
see how this arises, we first show (in section 3) that, if $f \in C_s(S^2)$, then the function
$e^{is\phi}f_I(\theta,\phi)$ extends continuously from $U_I$ to  $U_I \cup {\bf N}$,
and the function $e^{-is\phi}f_I(\theta,\phi)$ extends continuously from $U_I$ to
$U_I \cup {\bf S}$.  (This is shown by considering $f_R$ for those $R$ with ${\bf N}, {\bf S}
\in U_R$.)  Thus there is a well defined linear functional $L: C_s(S^2) \to \CC$, given by
$Lf = \lim_{\theta \to 0^+} e^{is\phi}f_I(\theta,\phi)$.  It is evident from (\ref{sylmfrm})
that $L({}_sY_{lm})$ must be zero if $m \neq -s$.  This indicates that $_sY_{l,-s}$ could have a
special status; in fact, it is a multiple of the spin $s$ zonal harmonic, in the sense of the following
results, which we prove in section 5.\\

{\bf Lemmas on Zonal Harmonics}
{\em
Define $_sZ_l = (-1)^{s^+}\sqrt{(2l+1)/4\pi}\ _sY_{l,-s}$.  Then:\\
(a) $L(_sZ_l) = (2l+1)/4\pi$.\\
(b) For any $f \in {\cal H}_{ls}$, $Lf = \langle f,\ _sZ_l \rangle$.\\
Next, let $P_{ls}$ be the projection in $L^2(S^2)$ onto ${\cal H}_{ls}$, so that $P_{ls}$ has kernel
$K^{ls}$, where for any rotations $R, R'$, $K^{ls}_{R',R}(x,y) = \sum_{m=-l}^l {}_sY_{lmR'}(x)\overline{{}_sY_{lmR}(y)}$
for $x \in U_{R'}$, $y \in U_{R}$.  Then:\\
(c) If $x \in U_I \cap U_{R'}$ and ${\bf N} \in U_R$, then
${}_{s}Z_{lI}(x) =
e^{is(\psi_1-\psi_2)}K^{ls}_{R',R}(x,{\bf N})$, where
$\psi_1 = \psi_{xR'I}$, and
$\psi_2$ is the angle from $\rho_R({\bf N})$  to $\partial/\partial y$ at ${\bf N}$.\\
(d) If $x \in U_R$, then $K^{ls}_{R,R}(x,x) = (2l+1)/4\pi$.}\\

In section 6, we begin our study of spin wavelets.  Before explaining our construction,
we need to explain the ideas which are used in the spin $0$ case on the sphere,
and on more general manifolds.  The references here are \cite{narc1}, \cite{narc2},
\cite{gmcw}, \cite{gmfr}, \cite{gm3}.

A word about notation: in sections 3-5 we often use the variable ``$f$'' to denote a spin function on the
sphere.  In sections 6-8 the variable ``$f$'' will be reserved for another purpose.

Specifically, say $f \in {\cal S}(\RR^+)$, $f \neq 0$,
and $f(0) = 0$.  Let  $({\bf M},g)$ be a smooth compact
oriented Riemannian manifold, and let $\Delta$ be the Laplace-Beltrami operator on
$\bf M$ (for instance, the spherical Laplacian if $\bf M$ is $S^2$).  Let $K_t$
be the kernel of $f(t^2\Delta)$.  Then, as we shall explain, the functions
\begin{equation}
\label{manwavdf0}
w_{t,x}(y) = \overline{K}_t(x,y),
\end{equation}
if multiplied by appropriate weights,
can be used as wavelets on ${\bf M}$.  In case ${\bf M} = S^n$ and $f$ has compact support
away from the origin, we shall say these wavelets are ``needlet-type''; the theory of needlets
was developed in \cite{narc1}, \cite{narc2}.  (Actually, the definition of needlet
in \cite{narc1}, \cite{narc2} is slightly different from (\ref{manwavdf}), as we
shall explain.)  The theory was worked out in full generality in \cite{gmcw}, \cite{gmfr}, \cite{gm3}.

Specifically, one starts with the {\em Calder\'on formula}: if $c \in (0,\infty)$ is defined by
$c = \int_0^{\infty} |f(t)|^2 \frac{dt}{t}$,
then for all $u > 0$,

\begin{equation}
\label{cald}
\int_0^{\infty} |f(tu)|^2 \frac{dt}{t} = c < \infty.
\end{equation}

Discretizing (\ref{cald}), if $a > 1$ is sufficiently close to $1$,
one obtains a special form of
{\em Daubechies' condition}:
for all $u > 0$,
\begin{equation}
\label{daub}
0 < A_a \leq \sum_{j=-\infty}^{\infty} |f(a^{2j} u)|^2 \leq B_a < \infty,
\end{equation}
where
\begin{align}
\label{daubest}
A_a =& \frac{c}{2|\log a|} \left(1 - O(|(a-1)^2 (\log|a-1|)|\right), \\\label{daubest2} B_a=&
\frac{c}{2|\log a|}\left(1 + O(|(a-1)^2 (\log|a-1|)|)\right).
\end{align}
((\ref{daubest})  and (\ref{daubest2})  were proved in \cite{gm1}, Lemma
7.6.  In particular, $B_a/A_a$ converges nearly quadratically to $1$ as $a \rightarrow 1$.  For example,
Daubechies calculated that if $f(u) = ue^{-u}$ and $a=2^{1/3}$, then $B_a/A_a = 1.0000$  to four
significant digits.

Let $P$ be the projection in $L^2(S^2)$ onto the space of constant functions (the null
space of $\Delta$.)
Suppose that $a > 1$ is sufficiently close to $1$ that $A_a$ and $B_a$
are close (or that $a,f$ have been chosen in such a manner that $B_a = A_a$).
Then
for an appropriate discrete set  $\{x_{j,k}\}_{(j,k)\in \ZZ\times [1,N_j]}$ on $\bf M$,
and certain weights $\mu_{j,k}$, the collection of $\{\phi_{j,k}:= \mu_{j,k}w_{a^j,x_{j,k}}\}$
constitutes a wavelet  frame for $(I-P)L^2(S^2)$.  By this we  mean that there exist ``frame bounds' $0<A\leq B<\infty$,
such that  for any $F\in (I-P)L^2(S^2)$ the following holds:
 \begin{align}\label{frame}
 A\parallel F\parallel_2^2 \leq \sum_{j,k} \mid \langle F, \phi_{j,k}\rangle\mid^2 \leq B\parallel F\parallel_2^2.
 \end{align}
If the points $\{x_{j,k}\}$ are selected properly, $B$ can be made arbitrarily close to $B_a$ and
$A$ can be made arbitrarily close to $A_a$.  Because of the proximity of $A_a$ to $B_a$,
the frame is therefore a ``nearly tight frame".  In passing from one scale to
another, we use ``dilations'' in the sense that we are looking at the kernel of $f(t^2\Delta)$ for
different $t$, which corresponds to dilating the metric on the manifold.  On the sphere, the wavelets at each scale
$j$ are all rotates of each other, up to constant multiples.  (This is analogous to the usual situation
in which all the wavelets at each scale are translates of each other.)

Again, $w_{t,x}$ is given by (\ref{manwavdf0}).  On the sphere $S^2$, we have
\begin{equation}
\label{sphkt0}
K_t(x,y) = \sum_{l=0}^\infty \sum_{m=-l}^l f(t^2\lambda_l)  Y_{lm}(x) \overline{ Y_{lm}}(y).
\end{equation}

The most important cases to consider on the sphere are the case in which $f$ has compact support away from $0$ (the
``needlet-type'' case, of Narcowich, Petrushev and Ward), and the case in which $f(u) = u^r e^{-u}$ for some
integer $r \geq 1$ (the ``Mexican needlet'' case, as considered in \cite{gmcw}, \cite{gmfr}).  (Actually, in their
work, Narcowich, Petrushev and Ward used $l^2$ in place of $l(l+1)$ in (\ref{sphkt}), but this is a minor distinction.)

Needlet-type wavelets and Mexican needlets each have their own advantages.
Needlet-type wavelets have these advantages: for appropriate $f$, $a$, $\mu_{j,k}$ and $x_{j,k}$, needlet-type wavelets
are a tight frame on the sphere (i.e., $A=B$ in   (\ref{frame})), and the frame elements at non-adjacent scales
are orthogonal.   In fact, for appropriate needlet-type $f$, $a$, $\mu_{j,k}$ and $x_{j,k}$, one can arrange that
$A_a=B_a=A=B=1$ in the discussion above.  In that case, the wavelets are called needlets.

Mexican needlets have their own advantages.
In \cite{gmcw}, \cite{gmfr}, an approximate formula is written down for them which can be used directly on the sphere.
(This formula, which arises from computation of a Maclaurin series, has been checked numerically.
It is work in progress, expected to be completed soon, to estimate the remainder terms in this
Maclaurin series.  In \cite{gmcw}, \cite{gmfr} the formula is written down only for $r=1$, but it can be readily
generalized to general $r$.  The formula indicates that Mexican needlets are quite analogous to
the Mexican hat wavelet which is commonly used on the real line.)
Mexican needlets have Gaussian decay at each scale, and they
do not oscillate (for small $r$).  Thus they can be implemented directly on the sphere, which is
desirable if there is missing data (such as the ``sky cut'' of the CMB).

It would be worthwile to utilize both needlets and Mexican needlets in the analysis of CMB, and the
results should be compared.  In this article we focus on needlet-type situations.

A key property of this kind of wavelet is its {\em localization}.
Working in the general situation (general ${\bf M}$, general $f$), one has:\\
 for every pair of
$C^\infty$ differential operators $X$ (in $x$) and $Y$  (in $y$) on ${\bf M}$, and for every nonnegative integer
$N$, there
 exists $c:=c_{N,X,Y}$ such that for all $t>0$ and $x, y\in {\bf M}$
\begin{equation}
\label{diagloc0}
\left| XYK_t(x,y) \right|\leq c ~ \frac{ t^{-(n+I+J)}}{(d(x,y)/t)^N},
\end{equation}
 where $I:=\deg X$, $J:=\deg Y$, and $d(x,y)$ denotes the geodesic distance from $x$ to $y$.
  In the case where ${\bf M}$ is the sphere, $f$ has compact support away from $0$,
and if one replaces $\lambda_l=l(l+1)$ by $l^2$ in the formula for $K_t$, this was
shown by Narcowich, Petrushev and Ward, in \cite{narc1} and \cite{narc2}.  The general result
was shown in \cite{gmcw}.

The proof of (\ref{diagloc0}) in \cite{gmcw} would in fact go through without change if, instead of
$K_t$ being the kernel of $f(t^2\Delta)$, where $\Delta$ is the Laplace-Beltrami operator, one assumed
that $K_t$ was the kernel of $f(t^2D)$, where $D$ is any smooth positive second-order elliptic differential operator,
acting on $C^{\infty}({\bf M})$.  It is natural to conjecture then that something analogous
holds if $D$ is instead a smooth positive second-order elliptic differential operator, between sections
of line bundles on ${\bf M}$.  Since $\Delta_s$ may be interpreted as being such an operator, it is natural
to want to prove the following result.\\

{\bf Theorem \ref{locestls}}{\em

Let $K_t$ be the kernel of $f(t^2\Delta_s)$.  Then:

For every $R, R' \in SO(3)$, every pair of compact sets ${\cal F}_R \subseteq U_R$ and ${\cal F}_{R'} \subseteq U_{R'}$,
and every pair of $C^{\infty}$ differential operators $X$ (in $x$) on $U_{R'}$ and $Y$ (in $y$) on $U_{R}$,
and for every nonnegative integer $N$, there exists $c$ such that for all $t > 0$, all $x \in {\cal F}_{R'}$
and all $y \in {\cal F}_{R}$, we have
\begin{equation}
\label{diaglocls0}
\left| XYK_{t,R',R}(x,y) \right|\leq c ~ \frac{ t^{-(n+I+J)}}{(d(x,y)/t)^N},
\end{equation}
where $I:=\deg X$ and $J:=\deg Y$.}\\

In this article, we only give the proof of Theorem \ref{locestls} in the easier ``needlet-type'' case,
where $f$ has compact support away from the origin.  The general case will be established elsewhere. \\

Let us now explain what we mean by spin wavelets.  For $R \in SO(3)$, and $x \in U_R$, let us define
\begin{equation}
\label{needdf0}
w_{txR} = \sum_{l \geq |s|} \sum_m \overline{f}(t^2\lambda_{ls})\:\overline{{}_sY_{lmR}}(x)\:{}_sY_{lm}.
\end{equation}
Then $w_{txR} \in C_s^{\infty}(S^2)$.  Also note that $w_{txR,R'}(y) = \overline{K}_{t,R,R'}(x,y)$
(if $x \in U_R$ and $y \in U_{R'}$);
this generalizes the case $s=0$ of (\ref{manwavdf0}).  Moreover, if
$F \in L^2_s(S^2)$, then for $x \in U_R$
\begin{equation}
\label{betFdf0}
\beta_{t,F,x,R} := \langle F, w_{txR}\rangle = (f(t^2\Delta_s)F)_R(x) = (\beta_{t,F})_R(x)
\end{equation}
where we have set $\beta_{t,F} = (f(t^2\Delta_s)F)$.  We call $\beta_{t,F,x,R}$ a {\em spin wavelet
coefficient} of $F$.  If $F = \sum_{l \geq |s|}\sum_m a_{lm}\:{}_sY_{lm} \in L^2_s(S^2)$, then
\begin{equation}
\label{betRexp0}
\beta_{t,F,x,R} = \sum_{l \geq |s|} \sum_m f(t^2\lambda_{ls})a_{lm}\:{}_sY_{lmR}(x).
\end{equation}

Let us now explain how, in analogy to the case $s=0$, one can obtain a nearly tight frame
from spin wavelets.  Let $P = P_{|s|,s}$
be the projection onto the ${\cal H}_{|s|,s}$, the null space of $\Delta_{ls}$ (in $C^{\infty}_s$).
Suppose that $a > 1$ is sufficiently close to $1$ that $A_a$ and $B_a$ (of (\ref{daubest}) and
(\ref{daubest2})) are close (or that $a,f$ have been chosen in such a manner that $B_a = A_a$).
We then claim that for an appropriate discrete set  $\{x_{j,k}\}_{(j,k)\in \ZZ\times [1,N_j]}$ on $S^2$,
certain $R_{j,k}$ with $x_{j,k} \in U_{R_{j,k}}$,
and certain weights $\mu_{j,k}$, the collection of $\{\phi_{j,k}:= \mu_{j,k}w_{a^j, x_{j,k}, R_{j,k}}\}$
constitutes a nearly tight wavelet frame for $(I-P)L^2_s(S^2)$.
(Note that this property is independent of the choice of
$R_{j,k}$.)  In that case, we will have that for some $C$, if $F \in (I-P)L^2_s$, then, in $L^2$,
$F \sim C\sum_{j,k}\langle F,\phi_{j,k} \rangle \phi_{j,k}$ (note that this sum is {\em independent}
of the choice of $R_{j,k}$).  The precise statement is Theorem \ref{framainsp} below, where we
explain how the $x_{j,k}$ and $\mu_{j,k}$ are to be chosen.  Theorem \ref{framainsp} will be proved 
elsewhere \cite{gms};
in this article, it is only used in section 8.

One can certainly choose $a$, and $f$ with compact support away from the origin,
so that $A_a = B_a = 1$ (recall (\ref{daub})).  Such $f$ were used in \cite{narc1} and \cite{narc2}.
It is then natural to conjecture, that as in the spin $0$ case of \cite{narc1} and \cite{narc2},
one can choose the points $x_{j,k}$ so that the $\{\mu_{j,k}\varphi_{j,k}\}$ are an exactly tight frame
for $(I-P)L^2_s(S^2)$, for suitable weights $\mu_{j,k}$.  One should note, however, that even
if such explicit points and weights could be found, it might not be practical for physicists to
use them in CMB analysis, where measurements are invariably taken at the HEALPix points \cite{healpix}.
This, then, is a situation where it is very useful (and perhaps essential) to have the flexibility in
choice of the $x_{j,k}$ which Theorem \ref{framainsp} allows.\\

In section 7, we begin to examine how our theory can be applied to cosmology.  We look at
random spin $s$ fields $G$ on $S^2$.
We say that $G$ is {\em isotropic} (in law)
if for every $x_1,\ldots,x_N \in U_I$, the joint probability distribution of
$G^R_I(x_1),\ldots,G^R_I(x_N)$ is independent of $R \in SO(3)$.
Similarly, if $G, H$ are both random spin $s$ fields, we say that $G,H$ are
{\em jointly isotropic} (in law) if for every $x_1,\ldots,x_N,y_1,\ldots,y_M \in U_I$,
the joint probability distribution of
$G^R_I(x_1),\ldots,G^R_I(x_N), G^R_I(y_1),\ldots,G^R_I(y_M)$ is independent of $R \in SO(3)$.

Say $G,H$ are jointly isotropic, and denote their spin $s$ spherical
harmonic coefficients by $a_{lm} = \langle G^R,\:_sY_{lm} \rangle$,
$b_{lm} = \langle H^R,\:_sY_{lm} \rangle$; these are random
variables.  In Theorem \ref{schurthm} we show that
$E(a_{lm}\overline{b}_{l'm'}) = 0$ unless $l=l'$ and $m=m'$ (see
\cite{bmv07} for the scalar case). Moreover
$E(a_{lm}\overline{b}_{lm})$ does not depend on $m$; we will denote
it by $C_{l,G,H}$.

In Theorem \ref{spinuncor}, we generalize a key result in \cite{BKMP06-3}, to show that the
spin wavelet coefficients of $G,H$ are asymptotically uncorrelated, under mild hypotheses.
That is, for $x \in U_R$, we define the spin wavelet coefficients
$\beta_{G,t,x,R} = <G,w_{txR}>$, $\beta_{H,t,x,R} = <H,w_{txR}>$, and we show that
$|\mbox{Cor}(\beta_{G,t,x,R},\beta_{H,t,y,R})| \to 0$ as $t \to 0^+$.  (Again we assume that
we are in the ``needlet-type'' situation, where $f$ has compact support away from $0$.)\\

Finally, in section 8, we offer a brief glimpse of the consequences of the considerations
of section 7 for stochastic limits and cosmology.

Say $\epsilon > 0$.
We use a nearly tight frame, as provided by Theorem \ref{framainsp}.  Specifically, in that theorem, we
fix a real $f$ supported in $[1/a^2,a^2]$ for which the Daubechies sum
$\sum_{j=-\infty}^{\infty} f^2(a^{2j} u) = 1$ for all $u > 0$, so that $A_a = B_a = 1$.
We then produce the $\{\phi_{j,k}\} = \{\phi_{j,k,\epsilon}\}$ with frame bounds $B = 1+\epsilon$ and
$A = 1-\epsilon$.

Say $G = \sum A_{lm}\:_sY_{lm}$
is an isotropic random spin $s$ field, and let $C_l = C_{l,G,G}$ (notation as in Theorem
\ref{spinuncor}).   For $j \in \ZZ$, we also let ${\gamma}_j = \sum_{a^{-2} \leq a^{2j}\lambda_{ls} \leq a^{2}} C_l(2l+1)$;
this is a small quantity since\\
$\sum_j \gamma_j \leq 3\sum_l C_l(2l+1) = 3\mbox{Var}G$.  (In all sums here, $l$ satisfies the implicit
restriction $l \geq |s|$.)
We let $\beta_{jk} = \beta_{jk\epsilon} = \langle G,\phi_{j,k,\epsilon}\rangle$.
We focus on the quadratic statistics

\begin{equation}
\label{gammatdf0}
\tilde{\Gamma}_{j} = \sum_{k}\left| \beta_{jk}\right| ^{2} = \sum_{k}\left| \beta_{jk\epsilon}\right| ^{2},\:\:\:\:\:\:
\widehat{\Gamma}_{j} = \sum_{l,m}f^{2}(a^{2j}\lambda _{ls})|A_{lm}|^2
\end{equation}
When considering $\tilde{\Gamma}_{j}$, we can even let $\epsilon$ depend on $j$ here.
In Proposition \ref{egam}, we show that
\begin{equation}
\label{egam0}
E(|\widehat{\Gamma}_{j}- \tilde{\Gamma}_{j}|) \leq \epsilon \gamma_j
\end{equation}
for all $j$.

Since we can let $\epsilon << 1$ depend on $j$, by Tchebychev's
inequality, the probability distribution of $\tilde{\Gamma}_{j}$ is
a small perturbation of that of $\widehat{\Gamma}_{j}$.

We conclude, then, by briefly discussing the latter.
We assume now that $G$ is {\em Gaussian}, by which we mean that
$\{\Re G_I(x): x \in U_I\} \cup \{\Im G_I(x): x \in U_I\}$ is a Gaussian family.
We shall also assume that $G = \sum A_{lm}\:_sY_{lm}$ is {\em involutive}, by which we mean that
$\overline{A}_{lm} = A_{l,-m}$ for all $m$ (so that, in particular, $A_{l0}$ is real).
This assumption is very natural to make in cosmology.  Indeed, as we shall explain in
section 4, Newman and Penrose explained how to write any random spin $s$ field
$F$ as $F = F_{\e} + iF_{\m}$, where the ``electric'' and ``magnetic'' parts
$F_{\e}$ and $F_{\m}$ are involutive.  In cosmology one assumes that these fields are
jointly isotropic.

On the other hand, the assumption of Gaussianity is common to a
large part of the literature on random fields, and it is predominant
in the cosmological area, both within physical and mathematical
articles. Concerning the physical literature, Gaussianity is
motivated by the widely dominating inflationary model for the Big
Bang dynamics, which predicts the polarization random field is the
outcome of Gaussian fluctuations of a quantum mechanical origin.
From the mathematical point of view, Gaussianity is nearly a
mandatory assumption in the \emph{high-frequency asymptotics}
framework we shall entertain here (compare \cite{mar06}); indeed, as
a consequence of the characterization of isotropic random fields
given in (\cite{bmv07}), the Gaussianity of the spherical harmonic
coefficients is equivalent to their stochastic independence, which
in turns make limit theorems feasible. The extension to non-Gaussian
circumstances is certainly a crucial issue for further research; we
stress, however, that in the framework of high-resolution
asymptotics no rigorous results are known in the literature in
non-Gaussian circumstances, even in the standard case of scalar
random fields.

Under these assumptions, in Proposition \ref{centlim} we prove a
central limit theorem for the $\widehat{\Gamma }_{j}$, namely that
\begin{equation*}
\frac{\widehat{\Gamma }_{j}-E\widehat{\Gamma }_{j}}{\sqrt{Var\left\{
\widehat{\Gamma }_{j}\right\} }}\rightarrow _{d}N(0,1)\text{ as }%
j\rightarrow -\infty \text{ ,}
\end{equation*}%
$\rightarrow _{d}$ denoting as usual convergence in probability law, and
$N(0,1)$ denoting as usual the standard normal distribution.

Let us explain the significance of this result. As we explained
earlier, our main motivating rationale is provided by the
observations of a polarization random field $F$ in cosmological
satellite experiments. We shall then consider the realistic
situation where a physical model is adopted to provide an expected
value for the ``angular power spectrum'' $C_l,$ of the random spin
field $F,$ typically as a function $C_{l}=C_{l}(\Psi) $, of some
vector of physical parameters $\Psi .$ In the case of cosmological
polarization data, the function $C_{l}(.)$ is estimated by the
numerical solution of partial differential (thermodynamic) equations
representing the evolution of a matter/energy fluid, and the parameters $%
\Psi $ are related to constants of fundamental physics. Propositions
\ref{egam} and \ref{centlim} highlight the possibility of using spin
quantities, such as the spin wavelet coefficients $\beta_{jk}$ for
$F,$ to provide asymptotically valid
estimators for scalar quantities of physical interest such as $%
C_{l}$. For instance, by Propositions \ref{egam} and \ref{centlim},
it is immediate to construct a test for the assumption that the
angular power spectra take the form $C_{l}(\Psi ),$ for some
physically motivated set of parameters which takes a given value
$\Psi =\Psi _{0}.$ It is indeed sufficient to focus on
\begin{equation*}
S_{j}:=\frac{\widehat{\Gamma }_{j}-E_{\Psi _{0}}\widehat{\Gamma }%
_{j}}{\sqrt{Var_{\Psi _{0}}\left\{ \widehat{\Gamma }_{j}\right\} } }\text{,}
\end{equation*}%
where $E_{\Psi _{0}}\widehat{\Gamma }_{j},Var_{\Psi _{0}}\left\{
\widehat{\Gamma }_{j}\right\} $ are the expected mean and variance
evaluated under the assumption that the vector of parameters is set
equal to $\Psi _{0}.$ Propositions \ref{egam} and \ref{centlim}
predict that, for $-j$ sufficiently large and $\epsilon$
sufficiently small (depending on $j$), the probability of
exceeding any given threshold value will be nearly provided by%
\begin{equation*}
\Pr \left\{ S_{j}\geq z_{\alpha /2}\right\} \sim \alpha \text{ ,}
\end{equation*}%
where $0<\alpha <1$ and $z_{\alpha /2}$ is the well-known quantile
of the standard Gaussian distribution $\Phi (.),$ i.e. $z_{\alpha
/2}:=\Phi ^{-1}(1-\alpha /2).$ Values of $S_{j}$ above $z_{\alpha
/2}$ provide evidence at \emph{significance level} $\alpha $ against
the assumption that $\Psi _{0}$ represent the physical parameters in
the model that generated the observations.\\

We stress that the direct use of $\widehat{\Gamma }_{j}$ is only feasible
when observations from the whole sky are available.  In cosmology, one is hampered
by the presence of unobserved sky regions, for instance due to the 
foreground emission by the Milky Way and other astrophysical sources. The localization 
properties of needlet-type spin wavelets
in the real and harmonic domains are then clearly most valuable because, outside the masked
regions, wavelet coefficients in the high-frequency limit are essentially unaffected
by the missing observations. Moreover, the 
uncorrelation properties of needlet coefficients in the Gaussian case make a sound 
asymptotic statistical theory feasible, even in the presence of unobserved regions. 
Furthermore, these dual localization properties make possible the search for features 
and cosmological asymmetries, as done in the scalar case by \cite{cooray}. We refer to \cite{ghmkp}
for more discussion and details.  (The reader should be aware that in the notation of \cite{ghmkp}, the 
high-frequency limit corresponds to letting $j \to \infty$, not $-\infty$ as in the present article.)

\section{Spin $s$ line bundles}

In the next section, we shall review the groundbreaking work of Newman and Penrose
\cite{newpen} on spin weight $s$.  They were writing for physicists, so it is necessary to place
their observations in a rigorous mathematical context.  We do this in a manner which remains very close to
the spirit of the work of Newman and Penrose, and uses very little machinery; for another approach,
using more machinery, see \cite{eastod}.  We also contribute some new results, some of which will be
crucial for our work.

Newman and Penrose work on a two-dimensional surface.  Writing for physicists,
they say that a quantity $\eta$ has spin weight $s \in \ZZ$, provided that, whenever a tangent vector $m$
at a point transforms under coordinate change
by $m' = e^{i\psi} m$, the quantity $\eta$, at that point, transforms by
$\eta' = e^{is\psi} \eta$.  Our first task is to put this notion into acceptable mathematical language.


Let ${\bf M}$ be an oriented Riemann surface.
The complex structure on ${\bf M}$ determines a unique
metric (inner product on real tangent vectors at each point), up to conformal equivalence (that is, up to a
constant multiple at each point).  (Indeed, locally, such a metric is the pullback through a chart of the
standard metric on $\CC$; different choices of charts lead to conformally equivalent metrics, since the
Jacobian of a conformal mapping is always a constant multiple of an orthogonal transformation.)   We always
use a metric on ${\bf M}$ in this conformal class (that is, a metric which agrees with one determined
by the complex structure, up to a constant multiple at each point) without further comment.  Let us call
this conformal class the {\em intrinsic conformal class} on ${\bf M}$.

If $p \in {\bf M}$, as usual we
let ${\bf M}_p$ denote the real tangent space at $p$.  Since ${\bf M}$ is oriented, for each $v \in {\bf M}_p$
we may select $Jv \in {\bf M}_p$ at that point, such that $(v,Jv)$ is
an oriented orthogonal basis for ${\bf M}_p$.  Surely $Jv$ is unique, up to a positive
multiple.  Thus, if $v,w$ are nonzero tangent vectors at a point,
and $g$ is our metric, we may very naturally define the
{\em angle from $v$ to $w$} to be the angle $\psi \in (-\pi,\pi]$ determined by the two properties: \\
\ \\
(1) $g(v,w) = r \cos \psi$ for some $r > 0$; and\\
(2) $\sgn\: \psi = \sgn\: g(Jv,w)$.\\
\ \\
Now, say we have an atlas ${\cal U}$ on ${\bf M}$.  Suppose that for each chart
$U_{\alpha} \in {\cal U}$ we have a section $\rho_{\alpha}$ of
$TU_{\alpha}$ (the real tangent bundle of $U_{\alpha}$).  (This section will give, at each point of $U_{\alpha}$,
a ``reference direction''.) Set $\rho = \{\rho_{\alpha}\}_{U_{\alpha} \in {\cal U}}$.  Say also $s \in \ZZ$.
We then define ${\bf L}^s$, the spin $s$ line bundle associated to $({\bf M},{\cal U}, \rho)$,
as follows.  Say $U_{\alpha}, U_{\beta}$ are charts in ${\cal U}$, and define, for each $p \in U_{\alpha} \cap U_{\beta}$,
$\psi_{p \beta \alpha}$ to be the angle from $\rho_{\alpha}(p)$ to $\rho_{\beta}(p)$.  We then set
\[ \lambda_{\beta \alpha}(p) = e^{is\psi_{p \beta \alpha}}. \]
Then $\lambda_{\alpha \beta}\lambda_{\beta \alpha} = 1$ on $U_{\alpha} \cap U_{\beta}$, and if
$U_{\gamma}$ is a chart in ${\cal U}$ as well, then
$\lambda_{\alpha \gamma}\lambda_{\gamma \beta}\lambda_{\beta \alpha} = 1$ on
$U_{\alpha} \cap U_{\beta} \cap U_{\gamma}$.

By a standard argument (see, for instance, the argument on \cite{horscv}, pages 139-140),
the $\lambda_{\beta \alpha}$
may be used as transition functions to define a complex
(not necessarily holomorphic) line bundle ${\bf L}^s$; this bundle is unique
up to isomorphism.
If $\pi$ is the projection in ${\bf L}^s$ onto ${\bf M}$,
there are therefore diffeomorphisms
$\Phi_{\alpha}: \pi^{-1}(U_{\alpha}) \rightarrow U_{\alpha} \times \CC$, such that if
$p \in U_{\alpha} \cap U_{\beta}$,
and if
\[\Phi_{\alpha}(p,\eta) = (p,z),\]
then
\[\Phi_{\beta}(p,\eta) = (p,e^{is\psi_{p \beta \alpha}}z).\]
Let us say that these maps $\{\Phi_{\alpha}\}$ are {\em induced} by $({\cal U},\rho)$.

Note that, if $\Omega$ is a open subset of ${\bf M}$, we may naturally restrict
${\bf L}^s$ to become a line bundle ${\bf L}^s_{\Omega}$ over $\Omega$, by using the charts
$\{U_{\alpha} \cap \Omega: U \in {\cal U}\}$,
and the section $\rho^{\prime}_{\alpha} := \rho_{\alpha}$ on $U_{\alpha} \cap \Omega$.

Suppose now that ${\bf L}_1^s$, ${\bf L}_2^s$ are determined, respectively, by
$({\bf M}^1,{\cal U}^1, \rho^1)$ and $({\bf M}^2,{\cal U}^2, \rho^2)$.
Say $F: {\bf M}^1 \rightarrow {\bf M}^2$ is holomorphic,
and is a local diffeomorphism. Then $F$ naturally
gives rise to a smooth map $F_s: {\bf L}_1^s \rightarrow{\bf L}_2^s$ as follows.
Say ${(\cal U}^1, \rho^1)$ induces the maps $\{\Phi^1_{\alpha}\}$,
and ${(\cal U}^2, \rho^2)$ induces the maps $\{\Phi^2_{\beta}\}$.
Say $p \in U^1_{\alpha}$, a chart in ${\cal U}^1$, and
$q = F(p) \in U^2_{\beta}$, a chart in ${\cal U}^2$.  The map $F$
gives rise to a mapping $F_{*p}: {\bf M}^1_p \rightarrow {\bf M}^2_q$.
Let $\psi$ be the angle from $F_{*p}[\rho^1_{\alpha}(p)]$ to $\rho^2_{\beta}(q)$. Say
now $(p,\eta) \in \pi^{-1}(U^1_{\alpha})$, and that
\[\Phi^1_{\alpha}(p,\eta) = (p,z).\]
Then we define $F_s(p,\eta) = (q,\zeta)$, where
\[\Phi^2_{\beta}(q,\zeta) = (q,e^{is\psi}z).\]
Since $F_*$ preserves angles, it is easy to see that this definition
is independent of choice of $\alpha, \beta$.  Further, if $I: {\bf
M}_1 \rightarrow {\bf M}_1$ is the identity map, so is $I_s: {\bf
L}_1^s \rightarrow {\bf L}_2^s$.

Moreover, say $F$ is actually a biholomorphic map.  Then, since
$F_*$ and $(F^{-1})_*$ preserve angles, one easily sees that
$(F_s)^{-1} = (F^{-1})_s$.  If $f$ is a (continuous) section of
${\bf L}_2^s$, we may naturally define a section $f^F$ of ${\bf
L}_1^s$ by
\begin{equation}
\label{indsecdf}
f^F = F^{-1}_s \circ f \circ F.
\end{equation}
If $G: {\bf M}^0 \rightarrow {\bf M}^1$ is another biholomorphic map, we then have
\begin{equation}
\label{circsfnctr}
(f^F)^G = f^{F \circ G}.
\end{equation}
\ \\
This is all quite abstract, but in fact we shall
study spin $s$ bundles only when the manifold is an open subset of the sphere $S^2$ or the complex
plane $\CC$.  On $S^2$, by choosing the atlas and $\rho$ well, we can impose a
very interesting and very explicit structure.

On $S^2$, realized as $\{(x,y,z) \in \RR^3: x^2+y^2+z^2 = 1\}$, we let ${\bf N}$
be the north pole $(0,0,1)$, and let ${\bf S}$ be the south pole, $(0,0,-1)$.
We define the chart $U_I$ to be $S^2 \setminus \{{\bf N}, {\bf S}\}$.  (The choice
of chart map from $U_I$ to $\CC$ is irrelevant; one can take it to be a stereographic
projection.)  We obtain all other charts in our atlas by rotating $U_I$.  Thus, if $R \in SO(3)$, we
define $U_R = RU_I$.  On $U_I$ we often use standard spherical coordinates $(\theta,\phi)$,
and at each point $p$ of $U_I$ we let $\rho_I(p)$ be the unit tangent vector at $p$
which is tangent to the circle $\theta = $ constant. (For definiteness, if this circle
is $x = r\cos \phi$, $y = r\sin \phi$, $z = $ constant, let us choose $\rho_I(p)$ to
point in the direction in which $\phi$ increases, the ``counterclockwise'' direction.)
On any other chart $U_R$, we choose
\begin{equation}
\label{rotsec}
\rho_R(Rp) = R_{*p} [\rho_{I}(p)].
\end{equation}
for any $p \in U_0$.

On the chart $U_R$ we can use coordinates $(\theta_R,\phi_R)$, obtained by rotation
of the $(\theta,\phi)$ coordinates.  Thus, if the
$(\theta,\phi)$ coordinates of $p \in U_I$ are $(\theta_0,\phi_0)$, then the
$(\theta_R,\phi_R)$ coordinates of $Rp \in U_R$ are also $(\theta_0,\phi_0)$.
Then, at any $q \in U_R$, $\rho_R(q)$ is the unit tangent vector to the circle
$\phi_R = $ constant, pointing in the direction of increasing $\phi_R$.
Again, the angle $\psi_{p R_2 R_1}$ is the angle from $\rho_{R_1}(p)$ to
$\rho_{R_2}(p)$.  We would clearly obtain the same angle if $\rho_R(p)$ had been
chosen as the unit vector pointing in the direction of increasing $\theta_R$, or in the direction of
any fixed linear combination of the unit vectors pointing in the directions of increasing
$\phi_R$ and increasing $\theta_R$.  In this precise sense, $\psi_{p R_2 R_1}$ effectively measures the
angle by which the tangent plane at $p$ is rotated if one uses the $(\theta_{R_2},\phi_{R_2})$
coordinates instead of the $(\theta_{R_1},\phi_{R_1})$ coordinates.   This is why
one should think of our choice of  $\rho_R(p)$ only as a convenient ``reference direction''.

On $\CC$ we will use only a single chart, $\CC$ itself, and if $p \in \CC$,
we let our reference direction $\rho(p)$ be $\partial/\partial y$.  Since there is only one chart, this
bundle is trivial.   However, if $F$ maps an open subset $U$ of $\CC$ conformally
onto another open subset $U$ of $\CC$, we can consider the map $F_s$, which, as we
shall see, can be quite interesting.\\

Now, if $f$ is a (continuous) section of ${\bf L}^s$ over
an open set $\Omega \subseteq S^2$, then
we may write
\[ \Phi_R(f(q)) = (q,f_R(q)) \]
for all $q \in U_R \cap \Omega$, for a suitable function $f_R: U_R \cap \Omega \rightarrow \CC$.
Note, then, that if
$p \in U_I \cap U_R \cap \Omega$, then
\begin{equation}
\label{rottrans1}
f_R(p) = e^{is\psi} f_I(p)
\end{equation}
where
\begin{equation}
\label{rottrans2}
\psi \mbox{ is the angle from } \rho_I(p) \mbox{ to } \rho_R(p)= R_{*p}(\rho_I(R^{-1}p)).
\end{equation}

Note also that if $\Omega \subseteq S^2$, and $R_0 \in SO(3)$,
any smooth function $f_{(R_0)}: U_{R_0} \cap \Omega
\rightarrow \CC$ determines a unique section $f$ of ${\bf L}^s$
over $U_{R_0} \cap \Omega$ with $f_{R_0} = f_{(R_0)}$.  Indeed,
say $R \in SO(3)$.  If $p \in U_{R_0} \cap U_R \cap \Omega$,
then $f_R(p)$ must be given by
\begin{equation}
\label{frrop}
f_R(p) = e^{is\psi}f_{R_0}(p).
\end{equation}
where of course, $f_{R_0} = f_{(R_0)}$, and where
\begin{equation}
\label{angrro}
\psi \mbox{ is the angle from } \rho_{R_0}(p) \mbox{ to } \rho_R(p).
\end{equation}
Conversely, there is surely a section $f$ of ${\bf L}^s$ over $U_{R_0} \cap \Omega$ with these $f_R$.\\

Consequently {\em the space of smooth sections} $f$ {\em of} ${\bf L}^s$ over $\Omega$ {\em may
be identified with the space }
\begin{equation}
\label{ident}
\{f = (f_R)_{R \in SO(3)}: \mbox{ each } f_R \in C^{\infty}(U_R \cap \Omega), \mbox{ and }
\mbox{ (\ref{frrop}), (\ref{angrro}) hold for all }
R_0, R \in SO(3) \mbox{ and all }
p \in U_{R_0} \cap U_R \cap \Omega\}.
\end{equation}
We shall frequently make this identification, which often greatly simplifies matters conceptually.
However, as is often the case in mathematics, certain properties are clearer, and more easily
verified, if one uses the coordinate-free line bundle point of view.

Equation (\ref{rottrans1}) describes what happens when one ``rotates the charts'', and
we now argue that ``you get the same answer if you rotate the section
instead''.
Thus, if $R \in SO(3)$, as we have seen, $R: S^2 \rightarrow S^2$ gives rise
to a smooth map $R_s: {\bf L}^s \rightarrow {\bf L}^s$.  If $f$ is a section
of ${\bf L}^s$ over $\Omega \subseteq S^2$, using (\ref{indsecdf}) we
obtain the ``rotated section''
\begin{equation}
\label{rotsecdf}
f^R = R_s^{-1} \circ f \circ R.
\end{equation}
over $R^{-1}\Omega$.  (Note that if $s=0$, $f^R = f \circ R$.) Thus, for $p \in \Omega$,
\[ (f^R)(R^{-1}p) = R_s^{-1} \circ f(p), \]
from which we see that
\begin{equation}
\label{rottrans3}
(f^R)_I (R^{-1}p) = e^{is\psi} f_I(p)
\end{equation}
where
\begin{equation}
\label{rottrans4}
\psi \mbox{ is the angle from } R_{*p}^{-1}(\rho_I(p)) \mbox{ to } \rho_I(R^{-1}p).
\end{equation}
Note that this $\psi$ is the same as the $\psi$ in (\ref{rottrans2}), since $R_{*p}$ preserves
angles.  Therefore
\begin{equation}
\label{rotsame}
f_R(p) = (f^R)_I(R^{-1}p).
\end{equation}
One could take the point of view of (\ref{ident}) here, and simply define
$f^R$ by setting $(f^R)_I(q) = f_R(Rq)$.  But then it would be tedious (though certainly possible) to verify that
the smoothness of $f$ implies the smoothness of $f^R$, or that $(f^{R_1})^{R_2} = f^{R_1R_2}$,
both of which are evident from (\ref{rotsecdf}).

We also need to discuss the metric and the orientation we will use on $S^2$.
Since rotations are holomorphic maps on $S^2$, any metric in the intrinsic conformal
class on $S^2$ is rotationally invariant, up to a constant multiple at each point.
Thus there is a truly rotationally invariant metric in this conformal class, obtained by
rotation of the metric at any particular point on the sphere.  We note that such a
rotationally invariant metric coincides (up to a constant multiple) with the
Euclidean metric, as restricted to tangent vectors on the sphere.  (That is,
for some positive constant $c$, if
$v=\sum_{j=1}^3 a_j \frac{\partial}{\partial x}_j$ and
$w=\sum_{j=1}^3 b_j \frac{\partial}{\partial x}_j$ are real tangent vectors
at any point on the sphere, their inner product there is $c <v,w>_E$, where
$<v,w>_E = \sum_{j=1}^3 a_jb_j$.)  Indeed, since $<v,w>_E$ is rotationally
invariant, this fact need only be checked at the South pole $\bf S$.  We can determine
a metric in the intrinsic conformal class on $S^2$, by pulling back the standard
metric on $\CC$ through stereographic projection from the North Pole onto the tangent plane to the
sphere at $\bf S$; denote this stereographic projection by $T$. It is geometrically evident,
and it is easily calculated (see the formula (\ref{sterprojE}) below for
$\sigma = T/2$), that $T_{*{\bf S}}(\partial/\partial x) = \partial/\partial x$,
and $T_{*{\bf S}}(\partial/\partial y) = \partial/\partial y$.
Thus the tangent vectors $\frac{\partial}{\partial x}$ and $\frac{\partial}{\partial y}$,
are, up to a choice of conformal factor at ${\bf S}$,
an orthonormal basis at ${\bf S}$, as desired.  We use $<v,w>_E$ as our metric
on $S^2$.  We choose the orientation on $S^2$ which makes $T$ orientation-preserving.
Thus $(\partial/\partial x, \partial/\partial y)$ is an oriented orthonormal basis at
${\bf S}$, and, consequently, $(\partial/\partial y, \partial/\partial x)$ is an
oriented orthonormal basis at ${\bf N}$.\\

In our first new result about the spin $s$ bundle over $S^2$, we take a section $f$ of ${\bf L}^s$,
and examine the behavior of $f_I$ as we approach the north and south poles.

\begin{theorem}
\label{limns}
Let $f$ be a continuous section of ${\bf L}^s$ over $S^2$.  Then the function
$e^{is\phi}f_I(\theta,\phi)$ extends continuously from $U_I$ to  $U_I \cup {\bf N}$,
and the function $e^{-is\phi}f_I(\theta,\phi)$ extends continuously from $U_I$ to
$U_I \cup {\bf S}$.
\end{theorem}
{\bf Proof}  Let us show the second statement; the proof of the first statement is almost identical.
For any $\phi \in [0,2\pi)$, let $f_{I,\phi}(\theta) = f_I(\theta,\phi)$.
It will suffice to show that $\lim_{\theta \rightarrow \pi^-} e^{-is\phi} f_{I,\phi}(\theta)$
exists, uniformly in $\phi$, and that this limit is independent of $\phi$.

For $\phi \in [0,2\pi)$, let $\gamma_{\phi}$ be the great arc $\gamma_{\phi}(\theta) = (\theta,\phi)$
$(0 \leq \theta \leq \pi)$.  If $p = \gamma_{\phi}(\theta)$ is on this great arc, then $\rho_I(p)$ points
in the direction $\gamma_{\phi}'(\theta) \times p$ (cross product); thus
\[ \lim_{\theta \rightarrow \pi^-} \rho_I(\gamma_{\phi}(\theta)) := v_{\phi} \]
exists, and this limit is uniform in $\phi$ (indeed, a rotation about ${\bf N}$ takes any one of these
limiting situations into any other).  Moreover, evidently, $v_{\phi}$ is a tangent vector at ${\bf S}$,
and given our choice of orientation at ${\bf S}$,
\begin{equation}
\label{v0vp}
\mbox{the angle from } v_0 \mbox{ to } v_{\phi} \mbox{ is } \phi.
\end{equation}

Now, let $R$ be any rotation of $S^2$ which takes ${\bf S}$ to a point other than ${\bf N}$ or ${\bf S}$.
We then have that, for $(\theta,\phi) \in U_I \cap U_R$,
\[ f_{I,\phi}(\theta) = e^{-is\psi_{\theta,\phi}}f_R(\theta,\phi) \]
where $\psi_{\theta,\phi}$ is the angle from $\rho_I(p)$ to $\rho_R(p)$, if $p = (\theta,\phi)$.
For fixed $\phi$, the limit
as $\theta \rightarrow \pi^-$ of the right side exists (uniformly in $\phi$) and equals
$e^{-is\Psi{\phi}}f_R({\bf S})$, where now
\begin{equation}
\label{psiphidf}
\Psi_{\phi} \mbox{ is the angle from } v_{\phi} \mbox{ to } \rho_R({\bf S}).
\end{equation}
Thus for any $\phi$,
\[
 \lim_{\theta \rightarrow \pi^-} e^{-is\phi} f_{I,\phi}(\theta)
= e^{-is(\phi + \Psi{\phi})}f_R({\bf S}) = e^{-is(\phi + \Psi{\phi}-\Psi_0)} e^{-is\Psi_0} f_R({\bf S})
= e^{-is(\phi + \Psi{\phi}-\Psi_0)} \lim_{\theta \rightarrow \pi^{-}} f_{I,0}(\theta)
= \lim_{\theta \rightarrow \pi^{-}} f_{I,0}(\theta),
\]
since $\Psi_{\phi}-\Psi_0$ equals $-\phi$ (mod $2\pi$), by (\ref{v0vp}) and (\ref{psiphidf}).
This proves the second statement; the proof of the first statement is entirely similar,
if one takes into account the difference in orientation at ${\bf N}$.

\section{The Newman-Penrose Theory}

In this section we review the spin $s$ theory of Newman and Penrose \cite{newpen}.  We do so
within the rigorous framework we have presented in section 3.  Most of our arguments are very close to theirs.
We resume presenting new results in the next section.

Although we earlier used stereographic projection $T$ from the north pole onto the tangent plane at ${\bf S}$, in
what follows we will follow Newman and Penrose and instead use stereographic
projection $\sigma$ from the north pole onto the equatorial plane, from $U_I \subseteq S^2$ to $\CC^*
:= \CC \setminus \{0\}$.  Note that $\sigma = T/2$ is still orientation-preserving.
Explicitly, if $p \in U_I$ has Euclidean coordinates $(x,y,z)$, then
\begin{equation}
\label{sterprojE}
\sigma p = \frac{1}{1-z}(x+iy),
\end{equation}
which one easily visualizes by using similar triangles.
From this, and the formulas $\sin \theta = 2 \sin\frac{\theta}{2}\cos\frac{\theta}{2}$,
$1-\cos\theta = 2\sin^2\frac{\theta}{2}$, it follows that
if $p \in U_I$ has standard spherical coordinates
$(\theta,\phi)$, then
\begin{equation}
\label{sterproj}
\sigma p = e^{i\phi}\cot(\frac{\theta}{2}).
\end{equation}

The map $\sigma$, being a conformal map, preserves angles; we will also need
to know the factor by which it multiplies lengths.  Thus, if $p \in U_I$,
since $\sigma$ is conformal, there is a number $\lambda_p > 0$ such that
for any tangent vector $v$ at $p$, $|\sigma_{*p}v| = \lambda_p <v,v>_E^{1/2}$.
We need to find the conformal factor $\lambda_p$.  For this, we may assume that $v = \rho_I(p)$ is
a unit tangent vector tangent to the circle $\phi = $ constant through $p$.
Restricted to that circle, by (\ref{sterprojE}), the map $\sigma$ is simply a
dilation by the factor $(1-z)^{-1}$, so we must have $\lambda_p = (1-z)^{-1}$.
It is easy to calculate, then, from (\ref{sterprojE}), that $\lambda_p =
\frac{1}{2}(1+|\sigma p|^2)$.  For this reason, for $\zeta \in \CC$, we set
\begin{equation}
\label{pzetdf}
P(\zeta) = \frac{1}{2}(1+|\zeta|^2),
\end{equation}
and we have shown that:
\begin{equation}
\label{conffacster}
\mbox{if } p \in U_I,\: \sigma p = \zeta, \mbox{ and } v \mbox{ is a tangent vector at }
p, \mbox{ then } |\sigma_{*p}v| = P(\zeta)<v,v>_E^{1/2}.
\end{equation}
We also remark the following fact (which is also geometrically evident).
Since $\sigma$ restricted to the aforementioned circle is just a dilation,
if
\begin{equation}
\label{rotsig}
\mbox{If } p = (\theta,\phi) \in U_I, \mbox{ then the angle from } \sigma_*\rho_I(p) \mbox{ to }
\partial/\partial y \mbox{ is } -\phi.
\end{equation}
Indeed, it equals the angle from $\rho_I(p)$ to $\partial/\partial y$ in $\RR^3$, which is surely
$-\phi$.\\

We are now ready for the observations of Newman and Penrose, expressed in our language.
These observations will be preceded by bullet points
in what follows.  We will provide rigorous proofs of their observations, within the context of spin $s$
line bundles.  Our proofs are modified versions of their arguments.\\

Say $\Omega \subseteq U_I$ is open.  Suppose $h_{(I)}: \Omega \rightarrow \CC$ is smooth.
Following Newman and Penrose we define
\[ \edth_{sI} h_{(I)} = -(\sin \theta)^s \left(\frac{\partial}{\partial \theta} +
\frac{i}{\sin \theta}\frac{\partial}{\partial \phi} \right)(\sin \theta)^{-s} h_{(I)}. \]
(Here we are again using standard spherical coordinates $(\theta,\phi)$ on $U_I$.)
Thus $\edth_{sI} h_{(I)}$ is a smooth function on $\Omega$.

Next, if $h$ is a smooth section of ${\bf L}^s_{\Omega}$, we define $\edth h$ to be the section of
${\bf L}^{s+1}$ over $\Omega$ such that
\[ (\edth h)_I = \edth_{sI} h_I. \]
Also, if $\Omega_0 \subseteq \CC^*$ is open, and $u$ is a smooth section of
${\bf L}^s_{\Omega_0}$, we define the smooth section ${\cal D}_s\:u$ of ${\bf L}^{s+1}_{\Omega_0}$
by
\[ {\cal D}_s\: u =  2P^{1-s}\frac{\partial P^s u}{\partial \zeta}.\]
(Here $P$ is as in (\ref{pzetdf}).
Recall that we have only one chart on $\CC^*$, so sections of
${\bf L}^s_{\Omega_0}$ or ${\bf L}^{s+1}_{\Omega_0}$ may be naturally identified with smooth
functions from $\Omega_0$ to $\CC$.)

If $f$ is a section of ${\bf L}^s_{\Omega}$, where $\Omega \subseteq U_I$, let us write $f^{\sigma^{-1}}) =
f \circ_s \sigma^{-1}$.  The first key observation of Newman and Penrose is:\\
\ \\
$\bullet$ Suppose $f$ is a smooth section of ${\bf L}^s_{\Omega}$, where $\Omega \subseteq U_I$ is
open.  Then
\begin{equation}
\label{edthgoP}
(\edth f)^{\sigma^{-1}} = {\cal D}_s\:(f^{\sigma^{-1}}).
\end{equation}

To see this, one notes that $\edth_{sI} f_I = s\cot \theta f_I
- Lf_I$, where $L = (\partial/\partial \theta) + (i/\sin\theta)(\partial/\partial \phi)$.
Using (\ref{rotsig}), we see that the left side of (\ref{edthgoP}) is
$e^{-i(s+1)\phi}[(s\cot \theta)f_I \circ \sigma^{-1}-(Lf_I)\circ \sigma^{-1}]$.
Let $u = f_I \circ \sigma^{-1}$.
A brief calculation, using (\ref{sterproj}), the chain rule, and the fact that
$\cot(\theta/2)/\sin\theta = \frac{1}{2}\csc^2(\theta/2)$, shows that $(Lf_I) \circ \sigma^{-1}
= -e^{i\phi}\csc^2(\theta/2)\partial u/\partial \zeta$.  On the other hand, if $p = (\theta,\phi) \in U_I$
and $\zeta = \sigma p$, then $P(\zeta) = \frac{1}{2}[1+\cot^2(\theta/2)] = \frac{1}{2}\csc^2(\theta/2)$.
Thus the left side of (\ref{edthgoP}) equals
\begin{equation}
\label{escpze}
e^{-i(s+1)\phi}\left[(s\cot \theta)u+2Pe^{i\phi}\frac{\partial u}{\partial \zeta}\right].
\end{equation}
On the other hand, the
right side of (\ref{edthgoP}) is
\begin{equation}
\label{psphuz}
2P^{1-s}\frac{\partial e^{-is\phi}P^s u}{\partial \zeta},
\end{equation}
which equals
$2P^{1-s}[\partial e^{-is\phi}P^s/\partial \zeta]u + 2Pe^{-is\phi}\partial u/\partial \zeta$.
We are left with showing that
\begin{equation}
\label{edthgoPout}
2P^{1-s}\frac{\partial e^{-is\phi}P^s}{\partial \zeta} = e^{-i(s+1)\phi}s\cot \theta.
\end{equation}
Surely $\partial P/\partial \zeta = \frac{1}{2}\overline{\zeta}$.  On the other hand,
since $e^{2i\phi}=\zeta/\overline{\zeta}$, we find
$2ie^{2i\phi}\partial \phi/\partial \zeta = 1/\overline{\zeta}$, so that
$\partial \phi/\partial \zeta = 1/(2i\zeta)$.  Thus the left side of (\ref{edthgoPout})
equals
\[ e^{-is\phi}[-sP/\zeta+s\overline{\zeta}] = se^{-is\phi}[(|\zeta|^2-P)/\zeta]
= se^{-i(s+1)\phi}[(|\zeta|^2-P)/\cot(\theta/2)], \]
so all that is left to show is that $(|\zeta|^2-P)/\cot(\theta/2) = \cot \theta$.
But\\ $(|\zeta|^2-P)/\cot(\theta/2) = [\cot^2(\theta/2)-1]/[2\cot(\theta/2)] = \cot \theta$
as desired.\\

Next, Newman and Penrose observe:\\
\ \\
$\bullet$ Say $R \in SO(3)$.  Suppose $f$ is a smooth section of ${\bf L}^s_{\Omega}$, where
$\Omega \subseteq U_I \cap U_R$ is
open.  Say $R \in SO(3)$.  Then
\begin{equation}
\label{edthgoR}
\edth (f^R) = (\edth f)^R.
\end{equation}
(Note that these sections are defined over $R^{-1}\Omega \subseteq U_I$.)

To see (\ref{edthgoR}), it suffices to show that $[\edth (f^R)]^{\sigma^{-1}} =
[(\edth f)^R]^{\sigma^{-1}}$.   (Note that these sections are defined over
$\Omega_1 := \sigma(R^{-1}\Omega)$.)  Set $u = f^{\sigma^{-1}}$, defined over
$\Omega_0 := \sigma \Omega$, and let
\begin{equation}
\label{invdf}
{\cal I} = \sigma \circ R \circ \sigma^{-1}: \Omega_1 \rightarrow \Omega_0.
\end{equation}
Then by (\ref{edthgoP}),
\[ [\edth (f^R)]^{\sigma^{-1}} = {\cal D}_s\:(u^{\sigma \circ R \circ \sigma^{-1}})
= {\cal D}_s\:u^{\cal I}. \]
On the other hand,
\[  [(\edth f)^R]^{\sigma^{-1}}=
(\edth f)^{\sigma^{-1} \circ \sigma \circ  R \circ \sigma^{-1}}
= ({\cal D}_s\:u)^{\cal I}. \]
Thus, we need only show that if $u$ is a section of $L^s_{\Omega_0} \subseteq \CC^*$, then
\begin{equation}
\label{DgoI}
{\cal D}_s\:u^{\cal I} = ({\cal D}_s\:u)^{\cal I}.
\end{equation}
(Note that these sections are defined over $\Omega_1$.)

Now $\cal I$ is a conformal mapping.  At any point $\zeta$ of $\Omega_1$, if ${\cal I}'(\zeta)
= re^{i\psi}$ ($\psi \in [-\pi,\pi)$), then we know that
\begin{equation}
\label{invang}
\mbox{the angle from } {\cal I}_{*\zeta}(\partial/\partial y)  \mbox{ to } \partial/\partial y \mbox{ is }
-\psi,
\end{equation}
while for any real nonzero tangent vector $v$ at $\zeta$, $|{\cal I}_{*\zeta}v| = r|v|$.  By
(\ref{conffacster}) and (\ref{invdf}), we conclude that
\begin{equation}
\label{invdrv}
{\cal I}'(\zeta) = \frac{P({\cal I}(\zeta))}{P(\zeta)}e^{i\psi}.
\end{equation}
Accordingly
\begin{equation}
\label{invdrvcnj}
P(\zeta)e^{i\psi} = \frac{P({\cal I}(\zeta))}{\overline{{\cal I}'(\zeta)}}.
\end{equation}
Turning now to (\ref{DgoI}), we see by (\ref{invdrvcnj}) that
\[ ({\cal D}_s\:u^{\cal I})(\zeta) =
2P(\zeta)^{1-s}\frac{\partial \left[(P(\zeta)e^{i\psi})^s (u \circ {\cal I})(\zeta)\right]}
{\partial \zeta} =
\frac{2P(\zeta)}{P(\zeta)^s}
\frac{\partial \left[([P^s u] \circ {\cal I})(\zeta)/{\overline{{\cal I}'(\zeta)^s}}\right]}
{\partial \zeta}. \]
The key point is that $\overline{{\cal I}'(\zeta)^s}$ is {\em antiholomorphic}, so that
\[ ({\cal D}_s\:u^{\cal I})(\zeta) =
\frac{2P(\zeta)}{[P(\zeta)\overline{{\cal I}'(\zeta)]}^s} \frac{\partial([P^s u] \circ {\cal I})(\zeta)} {\partial \zeta} =
\frac{2P(\zeta)}{P({\cal I}(\zeta))^s e^{-is\psi}} \left[\frac{\partial(P^s u)} {\partial \zeta} \circ {\cal I}(\zeta)\right]
{\cal I}'(\zeta). \]
Finally, then, by (\ref{invdrv}),
\[ ({\cal D}_s\:u^{\cal I})(\zeta) = 2P({\cal I}(\zeta))^{1-s} e^{i(s+1)\psi}
\left[\frac{\partial(P^s u)} {\partial \zeta} \circ {\cal I}(\zeta)\right]
= ({\cal D}_s\:u)^{\cal I}(\zeta). \]
This proves (\ref{DgoI}), as desired.\\
\ \\
Next, say $\Omega \subseteq U_R$ is open, so that $R^{-1}\Omega \subseteq U_I$.  Say $h_{(R)}: \Omega \rightarrow \CC$
is smooth; then $h_{(R)} \circ R: R^{-1}\Omega \rightarrow \CC$ is smooth, so we can apply $\edth_{sI}$ to it and obtain another
smooth function on $R^{-1}\Omega$.  Now we define
\begin{equation}
\label{edthRdf}
\edth_{sR} h_{(R)} = [\edth_{sI} (h_{(R)} \circ R)] \circ R^{-1},
\end{equation}
which is, like $h_{(R)}$ itself, a smooth function on $\Omega$.  We now deduce:\\
\ \\
$\bullet$ Suppose $f$ is a smooth section of ${\bf L}^s$ over an open set $\Omega \subseteq U_I \cap U_R$.  Then
\begin{equation}
\label{edthsecgoR}
(\edth f)_R = \edth_{sR} f_R.
\end{equation}
To see this, say $p \in \Omega$.  Then by (\ref{rotsame}) and (\ref{edthgoR}),
\[ (\edth f)_R(p) = [(\edth f)^R]_I(R^{-1}p) = (\edth f^R)_I(R^{-1}p). \]
Thus, by the definition of $\edth_{sI}$ and (\ref{rotsame}) again,
\[ (\edth f)_R(p) = (\edth_{sI} f^R_I)(R^{-1}p) = (\edth_{sI} [f_R \circ R])(R^{-1}p) = (\edth_{sR} f_R)(p), \]
as desired.\\
\ \\
We can reformulate the last result more appealingly as follows.
On the chart $U_R$ we again use coordinates $(\theta_R,\phi_R)$, where if the
$(\theta,\phi)$ coordinates of $p \in U_I$ are $(\theta_0,\phi_0)$, then the
$(\theta_R,\phi_R)$ coordinates of $Rp \in U_R$ are also $(\theta_0,\phi_0)$.
Then, in these coordinates on $U_R$,
\begin{equation}
\label{edthRxp}
\edth_{sR} f_R = -(\sin \theta_R)^s \left(\frac{\partial}{\partial \theta_R} +
\frac{i}{\sin \theta_R}\frac{\partial}{\partial \phi_R} \right)(\sin \theta_R)^{-s} f_R.
\end{equation}

So far we have only defined $\edth$ over the chart $U_I$, but (\ref{edthsecgoR}) enables
us to define $\edth f$ for any smooth section $f$ of ${\bf L}^s_{\Omega}$, where now $\Omega$ is any
open subset of $S^2$.  We simply let $\edth f$ be the smooth section of ${\bf L}^{s+1}_{\Omega}$, with
the property that if $p \in \Omega$ is in $U_R$, then
\begin{equation}
\label{edthdfall}
(\edth f)_R(p) = \edth_{sR} f_R(p).
\end{equation}
Note that this definition makes sense.  Indeed, over $\Omega \cap U_I$, this section agrees with the
section $\edth f$ we have previously defined.  If, instead, $p$ is the north or south pole, choose a
particular $R \in SO(3)$ with $p \in U_R$, and let $h$ be the section of ${\bf L}^s$
over $\Omega \cap U_R$ such that $h_R = \edth_{sR} f_R$.  Then on $\Omega \cap U_R \cap U_I$, $h$ agrees with $\edth f$
as previously defined.  Thus, for any other rotation $R'$ such that $p \in U_{R'}$,
$h_{R'}(q) = \edth_{sR'} f_{R'}(q)$ for all $q \in \Omega \cap U_R \cap U_{R'} \cap U_I$.  Letting $q \rightarrow p$,
we see that $h_{R'}(p) = \edth_{sR'} f_{R'}(p)$ as well.  This shows that the definition (\ref{edthdfall})
makes sense.

As a consequence we have that:
\begin{equation}
\label{edthglb}
\edth \mbox{ maps smooth sections of } {\bf L}^s \mbox{ to smooth sections of } {\bf L}^{s+1}.
\end{equation}
Here the sections could be defined over the entire sphere;
$\edth$ is consequently called the ``spin-raising operator''.

Note also that, in (\ref{edthgoR}), we may now relax the hypothesis that $\Omega \subseteq U_I \cap U_R$,
and only assume that $\Omega \subseteq S^2$ is open.  Indeed, both sides of (\ref{edthgoR})
agree on $\Omega \cap U_I \cap U_R$; since both sides are smooth sections of ${\bf L}^s_{\Omega}$, they
must agree on all of $\Omega$.

As is customary
for sections of line bundles, we let $C^{\infty}({\bf L}^s)$
denote the space of smooth sections of ${\bf L}^s$, defined over all of $S^2$.
(In section 2, we called this space $C^{\infty}_s(S^2)$.)
Newman and Penrose make the following important observation:
\begin{equation}
\label{edthinj}
\hspace{-1.5 in}
\bullet \mbox{ If } s < 0,
\mbox{ then } \edth: C^{\infty}({\bf L}^s) \rightarrow C^{\infty}({\bf L}^{s+1}) \mbox{ is injective.}
\hspace{1.5 in}
\end{equation}

To see this, say $f \in C^{\infty}({\bf L}^s)$, and that $\edth f = 0$.
By Theorem \ref{limns}, $|f_I|$ is bounded on $U_I$.
Let $u = f_I \circ \sigma^{-1}$; then $u$ is bounded on $\CC^*$.
As in (\ref{edthgoP}) and (\ref{psphuz}), we see that
\[ 0 = (\edth f)^{\sigma^{-1}} = {\cal D}_s\:(f^{\sigma^{-1}})
= 2P^{1-s}[\partial e^{-is\phi}P^s u/\partial \zeta]. \]
The function $v := e^{is\phi}P^s\overline{u}$ is therefore holomorphic on $\CC^*$.
It is bounded near $0$, so it extends holomorphically to $\CC$.  But, since
$s < 0$, $[P(\zeta)]^s \rightarrow 0$ as $\zeta \rightarrow \infty$, $v \equiv 0$
by Liouville's theorem.  Consequently $u$, and hence $f$, are both zero, as desired.\\
\ \\
Say now that $\Omega \subseteq S^2$ is open.  There is evidently a well-defined map
of {\em complex conjugation} from $C^{\infty}({\bf L}^s_{\Omega})$ to
$C^{\infty}({\bf L}^{-s}_{\Omega})$, taking a section $f$ in the former space to
a section $\overline{f}$ in the latter space, defined by
\[  (\overline{f})_R(p) = \overline{f_R(p)} \]
for all $p \in U_R \cap \Omega$ (for any $R \in SO(3)$).  We may then define
the ``spin-lowering operator'' $\bedth$ by
\begin{equation}
\label{bedthdf}
\bedth f = \overline{\edth \overline{f}},
\end{equation}
so that $\bedth: C^{\infty}({\bf L}^s_{\Omega}) \rightarrow C^{\infty}({\bf L}^{s-1}_{\Omega})$.
In the coordinates $(\theta_R,\phi_R)$ on $U_R$, we evidently have that
\[ (\bedth f)_R = \bedth_{sR} f_R, \]
where
\begin{equation}
\label{bedthRxp}
\bedth_{sR} f_R = -(\sin \theta_R)^{-s} \left(\frac{\partial}{\partial \theta_R} -
\frac{i}{\sin \theta_R}\frac{\partial}{\partial \phi_R} \right)(\sin \theta_R)^{s} f_R.
\end{equation}

If $\Omega_0 \subseteq \CC^*$ is open, and $u$ is a smooth section of
${\bf L}^s_{\Omega_0}$, we define the smooth section $\overline{\cal D}_s\:u$ of ${\bf L}^{s-1}_{\Omega_0}$
by $\overline{\cal D}_s\:u = \overline{{\cal D}_s\:\overline{u}}$; thus
\[ \overline{\cal D}_s\: u =  2P^{1+s}\frac{\partial P^{-s} u}{\partial \overline{\zeta}}.\]
By taking complex conjugates in (\ref{edthgoP}), (\ref{edthgoR}), and (\ref{edthinj}), we find:\\
\ \\
$\bullet$ Suppose $f$ is a smooth section of ${\bf L}^s_{\Omega}$, where $\Omega \subseteq U_I$ is
open.  Then
\begin{equation}
\label{bedthgoP}
(\bedth f)^{\sigma^{-1}} = \overline{\cal D}_s\:(f^{\sigma^{-1}}).
\end{equation}

$\bullet$ Suppose $f$ is a smooth section of ${\bf L}^s_{\Omega}$, where $\Omega \subseteq S^2$ is
open.  Say $R \in SO(3)$.  Then
\begin{equation}
\label{bedthgoR}
\bedth (f^R) = \bedth f^R.
\end{equation}

\begin{equation}
\label{bedthinj}
\hspace{-1.5 in}
\bullet \mbox{ If } s > 0,
\mbox{ then } \bedth: C^{\infty}({\bf L}^s) \rightarrow C^{\infty}({\bf L}^{s-1}) \mbox{ is injective.}
\hspace{1.5 in}
\end{equation}

One can define $L^2$ sections of ${\bf L}^s$, (or on any other
line bundle over a compact oriented manifold, for that matter), and define an $L^2$ norm on them,
by using a finite atlas on the bundle and a partition of unity on the manifold.
Different choices of atlas and partition of unity lead to the same space of $L^2$ sections with an
equivalent $L^2$ norm.  For the bundle ${\bf L}^s$ there is a very natural choice of $L^2$ norm,
equivalent to the $L^2$ norm obtained through such a construction.  Note first that if $f$ is an $L^2$ section
of ${\bf L}^s$, then (a representative of) $f$ is defined almost everywhere.  Say $f$ is defined at
$p \in U_R \cap U_{R'}$.  We may then write $\Phi_R(f(p)) = (p,f_R(p))$, $\Phi_{R'}(f(p)) = (p,f_{R'}(p))$,
$f_{R'}(p) = e^{is\psi} f_R(p)$, where $\psi$ is the angle from $\rho_R(p)$ to $\rho_{R'}(p)$.  Thus $|f_R(p)|^2$
is independent of the choice of $R$, as long as $p \in U_R$.  Accordingly
\[ |f(p)|^2 := |f_R(p)|^2 \:\:\: (p \in U_R) \]
is well-defined a.e. on $S^2$.  We then define the $L^2$ norm of $f$ by
\[ \|f\|_{L^2}^2 = \int_{S^2} |f|^2 \]
where the integral is with respect to the usual surface measure on $S^2$.  We let $L^2({\bf L}^s)$
denote the space of $L^2$ sections of ${\bf L}^s$, with this norm.

We may now show:\\
\ \\
$\bullet$ The formal adjoint of $\edth: C^{\infty}({\bf L}^{s}) \rightarrow C^{\infty}({\bf L}^{s+1})$ is
$-\bedth$.\\
\ \\
For this, we need to show that, whenever $f \in C^{\infty}({\bf L}^{s})$ and $g \in C^{\infty}({\bf L}^{s+1})$,
we have
\begin{equation}
\label{edtadj}
<\edth f, g> = -<f,\bedth g>,
\end{equation}
where $<\cdot,\cdot>$ denotes $L^2$ inner product.

For this, we first note the following useful proposition:
\begin{proposition}
\label{s2notun}
(a) Let $N$ be a positive integer.  Then $S^2$ is not the union of $N$ open subsets of diameter less than
$\pi/N$.\\
(b) Any open set $\Omega \subseteq S^2$ of diameter less that $\pi/2$ is contained in some $U_R$.
\end{proposition}
{\bf Proof} For (a), say it were, and let $U_1$ be one of those open subsets.  By connectedness of $S^2$, $U_1$ must
intersect another one of those sets, say $U_2$.  Again by connectedness, $U_1 \cup U_2$ must intersect
another one of those sets, say $U_3$.  Continuing in this manner, we obtain the contradiction

\[ \pi = \diam S^2 \leq \diam U_1 + \diam U_2 + \ldots < N\pi/N. \]
For (b), we need to show that $\Omega$ omits some pair of antipodal
points in $S^2$.  If it did not, then if $A$ is the antipodal map, $\Omega \cup A\Omega = S^2$.
This contradicts (a).\\
\ \\
Let us then cover $S^2$ by a finite collection of open sets of diameter less than $\pi/5$,
choose a partition of unity $\{\zeta_j\}$ subordinate to this cover, and break up $f = \sum \zeta_j f$ and
$g = \sum \zeta_j g$.  By Proposition \ref{s2notun} (b), by using this breakup,
in proving (\ref{edtadj}), we may assume that $\supp f \cup \supp g \subseteq U_R$ for some $R$.  Then we need to show that
\[ \int_0^{2\pi} \int_{-\pi}^{\pi} [(\edth_{sR} f_R)\overline{g_R}](\theta_R,\phi_R)\sin \theta_R d\theta_R d\phi_R
= - \int_0^{2\pi} \int_{-\pi}^{\pi} [f_R\overline{\bedth_{s+1,R} g_R}](\theta_R,\phi_R)\sin \theta_R d\theta_R d\phi_R. \]
But this is an elementary calculation, using (\ref{edthRxp}) and (\ref{bedthRxp}).\\
\ \\
We are now ready (again following Newman and Penrose) to define the spin $s$ spherical harmonics.
We let $\{Y_{lm}: 0 \leq l, -l \leq m \leq l\}$ be the usual orthonormal basis of spherical harmonics
for $S^2$,
\begin{equation}
\label{sphareg}
Y_{lm}(\theta,\phi) = e^{im\phi}(-1)^{m^{-}}\:\frac{a_{lm}}{l!}
\sin^{2l}(\theta/2)\sum_r \left( \displaystyle^{l}_{r}\right)\left( \displaystyle^{\:\:\:l}_{r-m}\right)
(-1)^{l-r}\cot^{2r-m}(\theta/2).
\end{equation}
where
\[ a_{lm} = \left[(l+m)!(l-m)!(2l+1)/(4\pi)\right]^{\frac{1}{2}},\]
and our convention is that the binomial coefficient $\left(\displaystyle^n_j \right)$ is zero
if $j < 0$ (or $n < 0$).  (As usual $m^+ = \max(m,0)$, $m^{-} = max(-m,0)$.  The factor of
$(-1)^{m^{-}}$ in (\ref{sphareg}) is usually not included, but by including it, we have
$\overline{Y_{lm}} = Y_{l,-m}$, which will simplify many formulas in the sequel.)

For $l \geq |s|$, Newman and Penrose then define the spin $s$ spherical harmonic ${}_sY_{lm}$ by
\begin{equation}
\label{sylmdf}
{}_sY_{lm} =
\begin{cases}
\left[\frac{(l-s)!}{(l+s)!}\right]^{\frac{1}{2}}\edth^s Y_{lm}& (0 \leq s \leq l),\\
\left[\frac{(l+s)!}{(l-s)!}\right]^{\frac{1}{2}}(-\bedth)^{-s} Y_{lm} & (-l \leq s \leq 0).
\end{cases}
\end{equation}
The ${}_sY_{lm}$ are not defined for $l < |s|$.  Note that ${}_sY_{lm}$ is a smooth section of
${\bf L}^s$.  We have (\cite{newpen}, \cite{gmnrs}):
\begin{equation}
\label{sphaspn}
{}_sY_{lmI}(\theta,\phi) = e^{im\phi}(-1)^{m^{-}}\frac{a_{lm}}{[(l+s)!(l-s)!]^{1/2}}
\sin^{2l}(\theta/2)\sum_r \left( \displaystyle^{l-s}_{\:\:\:r}\right)\left( \displaystyle^{\:\:\:l+s}_{r+s-m}\right)
(-1)^{l-r-s}\cot^{2r+s-m}(\theta/2).
\end{equation}
To see this, let ${}_sv_{lm}$ denote the section of ${\bf L}^s$
over $U_I$ such that ${}_sv_{lmI}$ equals the right side of (\ref{sphaspn}).  We need to show that
${}_sv_{lm} = {}_sY_{lm}$ over $U_I$.  Note that this is surely true when $s = 0$.  Let
${}_su_{lm} = {}_sv_{lm}^{\sigma^{-1}}$.  By (\ref{edthgoP}), when $s \geq 0$, we need only show that
\begin{equation}
\label{ds0ulm}
{\cal D}^s\: {}_0u_{lm} = \left[\frac{(l+s)!}{(l-s)!}\right]^{\frac{1}{2}}{}_su_{lm}.
\end{equation}
(Here ${\cal D}^s = {\cal D}_{s-1}\: \circ \ldots \circ {\cal D}_0$.)
However, using (\ref{rotsig}) again, one easily calculates that, whenever $l \geq |s|$,
\begin{equation}
\label{sulmxp}
{}_su_{lm}(\zeta) = \frac{(-1)^{l-m^+}\: a_{lm}}{[(l+s)!(l-s)!]^{1/2}}(1+|\zeta|^2)^{-l}
\sum_r \left( \displaystyle^{l-s}_{\:\:\:r}\right)\left( \displaystyle^{\:\:\:l+s}_{r+s-m}\right)
\zeta^r(-\overline{\zeta})^{r+s-m}.
\end{equation}
From this, a brief calculation using Pascal's rule and the identities
\[ (l-s-r)\left( \displaystyle^{l-s}_{\:\:\:r}\right) = (r+1)\left( \displaystyle^{l-s}_{r+1}\right)
= (l-s)\left( \displaystyle^{l-s-1}_{\:\:\:\:\:\:r}\right) \]
shows that, whenever $l \geq |s|$,
\begin{equation}
\label{dsulm}
{\cal D}_s\:\:{}_su_{lm}(\zeta) = \left[(l-s)(l+s+1)\right]^{\frac{1}{2}}{}_{s+1}u_{lm}.
\end{equation}
(Here, our convention is, that since the factor $(l-s)$ appears on the right side of (\ref{dsulm}),
the right side is to be interpreted as zero if $l=s$.)  From this, (\ref{ds0ulm}) follows at once
by induction on $s$, for $s \geq 0$, and so we have (\ref{sphaspn}) for $s \geq 0$ as well.  We then
obtain (\ref{sphaspn}) for $s < 0$ by complex conjugation and use of the identity
\[ \overline{{}_sY_{lmI}} = \:{}_{-s}Y_{l,-m,I}. \]
Even better, by (\ref{dsulm}) and (\ref{edthgoP}), we have that, whenever $l \geq |s|$,
\begin{equation}
\label{edylm}
\edth \:{}_sY_{lm}(\zeta) = \left[(l-s)(l+s+1)\right]^{\frac{1}{2}}{}_{s+1}Y_{lm}.
\end{equation}
By taking complex conjugates of this, we also find that, whenever $l \geq |s|$,
\begin{equation}
\label{bedylm}
\bedth \:{}_sY_{lm}(\zeta) = -\left[(l+s)(l-s+1)\right]^{\frac{1}{2}}{}_{s-1}Y_{lm}.
\end{equation}
Newman and Penrose next observe that, for each $s$, the $\{{}_sY_{lm}\}$ are orthonormal.  Indeed,
set
\begin{equation}
\label{blsdf}
b_{ls} =
\begin{cases}
[(l+s)!/(l-s)!]^{1/2} & \mbox {if } 0 \leq s \leq l,\\
[(l-s)!/(l+s)!]^{1/2} & \mbox {if } -l \leq s \leq 0.
\end{cases}
\end{equation}
By (\ref{edylm}), if $s > 0$,
\begin{equation}
\label{ylmblsed}
{}_sY_{lm} = \frac{1}{b_{ls}}\edth^s Y_{lm},\:\:\:\:\:\: Y_{lm} = \frac{1}{b_{ls}}(-\bedth)^s\:{}_sY_{lm},
\end{equation}
while if $s < 0$,
\begin{equation}
\label{ylmblsbed}
{}_sY_{lm} = \frac{1}{b_{ls}}(-\bedth)^{-s} Y_{lm},\:\:\:\:\:\: Y_{lm} = \frac{1}{b_{ls}}\edth^{-s}\:{}_sY_{lm}.
\end{equation}

If $s \geq 0$, we see that
\[ <{}_sY_{lm},\:{}_sY_{l'm'}> = \frac{1}{b_{ls}}<\edth^s Y_{lm},\:{}_sY_{l'm'}>
= \frac{1}{b_{ls}}<Y_{lm},(-\bedth)^s{}_sY_{l'm'}> = <Y_{lm},Y_{l'm'}>= \delta_{ll'}\delta_{mm'}, \]
which proves the orthonormality.  If $s < 0$, one obtains the orthonormality through complex conjugation.

Thus, the $\{\:{}_sY_{lm}\}$ are an orthonormal set of eigenfuntions for the
(formally) positive self-adjoint operator $-\bedth \edth: C^{\infty}({\bf L}^s) \rightarrow
C^{\infty}({\bf L}^s)$.  Since this operator is elliptic (by (\ref{edthRxp}) and
(\ref{bedthRxp})), it has an orthonormal basis of smooth eigenfunctions, so at this point one
would surely expect:
\begin{equation}
\label{ylmob}
\hspace{-.05cm}
\bullet \mbox{ The } \{\:{}_sY_{lm}\} \mbox{ are an orthonormal basis of smooth eigenfunctions for }
-\bedth \edth: C^{\infty}({\bf L}^s) \rightarrow C^{\infty}({\bf L}^s).
\hspace{.05 cm}
\end{equation}
This would follow if we knew:
\begin{equation}
\label{edthsur}
\hspace{-1.67 in}
\bullet \mbox{ If } s > 0,
\mbox{ then } \edth^s: C^{\infty} \rightarrow C^{\infty}({\bf L}^{s}) \mbox{ is surjective.}
\hspace{1.67 in}
\end{equation}
Indeed, say $s > 0$.  Now $\edth^s: C^{\infty} \rightarrow C^{\infty}({\bf L}^{s})$ is continuous (using their
Fr\'echet space topologies).
(As usual, this topology is obtained by using a finite atlas on the bundle and a partition of unity
on the manifold.  Different choices of atlas and partition of unity lead to an equivalent Fr\'echet space
topology.)  Thus, since
the linear span of the $\{Y_{lm}\}$ is dense in $C^{\infty}$, the linear span of the
$\{\:{}_sY_{lm}\}$ is dense in $\edth^s C^{\infty} = C^{\infty}({\bf L}^{s})$, which is in turn
$L^2$-dense in the space of $L^2$ sections.  This would prove (\ref{ylmob}) for $s > 0$ (and hence for
all $s$ by complex conjugation), if we knew (\ref{edthsur}).  Moreover, (\ref{edthsur}) itself follows
from the injectivity of $\bedth^s: C^{\infty}({\bf L}^{s}) \rightarrow C^{\infty}$ (see
(\ref{bedthinj})), and by use of elliptic theory and Sobolev spaces.

Rather than proceed in this manner, and in order to extract important additional information,
we give a direct proof of (\ref{ylmob}) and (\ref{edthsur}), by
following the method of Newman and Penrose.  Their method was somewhat formal, but it can easily
be made rigorous.  In fact, we shall show:

\begin{theorem}
\label{cinfcgs}
Fix an integer $s$.\\
(a) Let $\{\|\:\|_N\}$ be a family of norms defining the Fr\'echet topology on $C^{\infty}({\bf L}^s)$.
Then for any $N$ there exist $C, M$ such that
\begin{equation}
\label{ylmpolygro}
\|{}_sY_{lm}\|_N \leq C(l+1)^M
\end{equation}
whenever $l \geq |s|$.\\
(b) For every $f \in C^{\infty}({\bf L}^{s})$, there exist numbers
$\{A_{lm}: l \geq |s|, -l \leq m \leq l\}$
such that
\begin{equation}
\label{cinfsxp}
f = \sum_{l \geq |s|}\sum_m A_{lm}\:{}_sY_{lm}
\end{equation}
(convergence in $C^{\infty}({\bf L}^{s})$), such that the $A_{lm}$ decay rapidly, by which we
mean that for every $N$ there exist $C, M$ such that
\begin{equation}
\label{almdec}
|A_{lm}| \leq C(l+1)^{-M}
\end{equation}
for all $l,m$.  Further, we must have $A_{lm} = <f,\:{}_sY_{lm}>$.
\end{theorem}
{\bf Proof} (a) is well-known when $s=0$.  It then follows for all $s$ by use of
(\ref{sylmdf}), as well as (\ref{bedthRxp}) and (\ref{edthRxp}) for general $R$.
Note, then, that any series as on the right side of (\ref{cinfsxp}), where the
$A_{lm}$ decay rapidly, does converge in $C^{\infty}({\bf L}^{s})$).

(b) is also well-known when $s=0$.  To prove it in general, we may as usual
assume $s > 0$.  In that case we may write
\[ \bedth^s f = \sum_{l \geq 0}\sum_m B_{lm}Y_{lm}
= \sum_{l \geq s}\sum_m B_{lm}Y_{lm} + g \]
where the $B_{lm}$ decay rapidly, and where $g = \sum_{l < s}\sum_m B_{lm}Y_{lm}
\in \mbox{ker}\edth^s$.  Now, with $b_{ls}$ as in (\ref{blsdf}), set $A_{lm}=(-1)^sB_{lm}/b_{ls}$
(for $l \geq |s|$).  Then, again using (a) for $s = 0$, we find
\begin{equation}
\label{bedfsmg}
\bedth^s (f - \sum_{l \geq s}\sum_m A_{lm}\:{}_sY_{lm}) = g \in \mbox{ker}\edth^s.
\end{equation}
Now it is a simple general fact, that if $D$ is a differential operator mapping smooth sections of
one hermitian line bundle, over a compact oriented manifold, to smooth sections of another
hermitian line bundle over that manifold,
and if $DF$ is in the kernel of $D^*$ for some smooth $F$, then $DF = 0$.  (Indeed, $D^*DF=0$,
so the inner product of $DF$ with itself is zero.)  Thus the left side of (\ref{bedfsmg}) is zero,
and so we have (\ref{cinfsxp}), by (\ref{bedthinj}).
Surely $A_{lm} = <f,\:{}_sY_{lm}>$, by the orthonormality of the ${}_sY_{lm}$.
This completes the proof.  Moreover, (\ref{edthsur}) follows as well, since if $s > 0$,
\[ f = \edth^s(\sum_{l \geq s}\sum_m C_{lm}Y_{lm}), \]
if $C_{lm} = A_{lm}/b_{ls}$.

As we have said, the completeness of the $\{{}_sY_{lm}\}$ as an orthonormal basis follows from
(\ref{edthsur}).  It also follows directly from Theorem \ref{cinfcgs} (b), which shows that the
linear span of the $\{\:{}_sY_{lm}\}$ is dense in $C^{\infty}({\bf L}^{s})$.

If $s > 0$, we may write

\begin{equation}
\label{cinfhg}
C^{\infty} = V \bigoplus W
\end{equation}
where
\begin{equation}
\label{vbedcn}
V = \bedth^s C^{\infty}({\bf L}^{s}) = \overline{ \langle Y_{lm}: l \geq s\rangle}
\end{equation}
and
\begin{equation}
\label{wkered}
W = \mbox{ker}\edth^s = \langle Y_{lm}: l < s\rangle.
\end{equation}
Moreover,
\begin{equation}
\label{edsbij}
\edth^s: V \rightarrow C^{\infty}({\bf L}^{s}) \mbox{ bijectively}.
\end{equation}
The statements (\ref{cinfhg}) and (\ref{edsbij}) are evident, once one identifies
$V$ and $W$ as spaces spanned by certain spherical harmonics, as in (\ref{vbedcn}) and
(\ref{wkered}).  One could obtain (\ref{cinfhg}) and (\ref{edsbij}) by using only the first
definitions of $V$ and $W$ in (\ref{vbedcn}) and (\ref{wkered}), by using elliptic theory
and (\ref{bedthinj}).  However, by identifying them as spaces spanned by certain spherical
harmonics, Newman and Penrose observe a crucial additional structure.

Namely, the decomposition (\ref{cinfhg}) {\em respects the splitting of a function into
its real and imaginary parts}.  That is, if $g \in V$ (or $W$), so are $\Re g$ and $\Im g$.
(This is evident, since $\overline{Y_{lm}}=Y_{l,-m}$.)
This splitting of any function in $V$ then induces, through the isomorphism in
(\ref{edsbij}), a corresponding splitting of any element of $C^{\infty}({\bf L}^s)$.
Thus, say $f \in C^{\infty}({\bf L}^{s})$.  Suppose
$s > 0$.  We write $f = \edth^s g$ for some $g \in C^{\infty}$, then
set
\[ f_{\e} = \edth^s \Re g,\:\:\: f_{\m} = \edth^s \Im g. \]
$f_{\e}$ and $f_{\m}$ are well-defined, independent of choice of $g$.
Indeed, if $\edth^s h = 0$, then $h \in W$, so
$\Re h$ and $\Im h$ are in $W$, so $\edth^s \Re h = \edth^s \Im h = 0$.
In particular, we could let $g$ be the unique solution of $f = \edth^s g$ in $V$.

Similarly, if $s < 0$, we write $f = \bedth^s g$ for some $g \in C^{\infty}$, then
set
\[ f_{\e} = \bedth^s \Re g,\:\:\: f_{\m} = \bedth^s \Im g. \]

Surely we always have $f = f_{\e} + if_{\m}$.

Again suppose $s > 0$.  By repeated use of (\ref{bedthgoR}), one sees that
the space $V$ is rotationally invariant, and then by repeated use of (\ref{edthgoR}),
one sees that
\begin{equation}
\label{emcircr}
(f^R)_{\e} = (f_{\e})^R,\:\:\: (f^R)_{\m} = (f_{\m})^R,
\end{equation}
for any $R \in SO(3)$.  Similarly these relations hold for $s < 0$ as well.

In fact, if $f$ is as in (\ref{cinfsxp}), it is very easy to find the expansions of
$f_{\e}$ and $f_{\m}$ in terms of the ${}_sY_{lm}$.  Say first $s > 0$.  The
solution $g$ of $f = \edth^s g$, with $g \in V$, is simply
\[ \sum_{l \geq s}\sum_m \frac{A_{lm}}{b_{ls}}Y_{lm}. \]
A very brief calculation, using $\overline{Y_{lm}} = Y_{l,-m}$, now shows that
\begin{equation}
\label{cinfsxpe}
f_{\e} = \sum_{l \geq |s|}\sum_m A_{lm\e}\:{}_sY_{lm}
\end{equation}
and
\begin{equation}
\label{cinfsxpm}
f_{\m} = \sum_{l \geq |s|}\sum_m A_{lm\M}\:{}_sY_{lm}
\end{equation}
where
\begin{equation}
\label{almexp}
A_{lm\e} = (A_{lm}+\overline{A_{l,-m}})/2,
\end{equation}
while
\begin{equation}
\label{almmxp}
A_{lm\M} = -i(A_{lm}-\overline{A_{l,-m}})/2.
\end{equation}
Similarly, (\ref{cinfsxpe}), (\ref{cinfsxpm}), (\ref{almexp}) and (\ref{almmxp}) hold
for $s < 0$ as well.
Up to normalization constants, these formulas agree with formulas in section II of
\cite{zs}, and also with equation (54) of \cite{cc} (to see the latter, use equations
(55), (65), (66) and (69) of that article).  (Those articles take $s=2$ and write $f=Q+iU$.
\cite{cc} uses the notation $a^G_{(lm)}$ and $a^C_{(lm)}$ in place of our $A_{lm\e}$ and $A_{lm\M}$.)\\
\ \\
Let us now let
\begin{equation}
\label{delsdf}
\Delta_s =
\begin{cases}
-\bedth\edth & \mbox{ if } s \geq 0,\\
-\edth\bedth & \mbox{ if } s < 0.
\end{cases}
\end{equation}
(Here $\Delta_s$ acts on smooth sections of ${\bf L}^s$.)
We also let
\begin{equation}
\label{evdf}
\lambda_{ls} =
\begin{cases}
(l-s)(l+s+1) & \mbox{ if } 0 \leq s \leq l,\\
(l+s)(l-s+1) & \mbox{ if } -l \leq s < 0.
\end{cases}
\end{equation}
We also let
\begin{equation}
\label{hlsdf}
{\cal H}_{ls} = <{}_sY_{lm}: -l \leq m \leq l> \mbox{  (for } l \geq |s|).
\end{equation}
Also
\begin{equation}
\label{plsdf}
P_{ls} \mbox{ denotes the projection in } L^2({\bf L}^s) \mbox{ onto } {\cal H}_{ls}.
\end{equation}
Then $L^2({\bf L}^s) = \bigoplus_{l \geq |s|} {\cal H}_{ls}$.  We also have that ${\cal H}_{ls}$
is the space of smooth sections of ${\bf L}^s$ which are eigenfunctions of
$\Delta_s$ for the eigenvalue $\lambda_{ls}$.  (This follows from (\ref{edylm}), (\ref{bedylm}), and
Theorem \ref{cinfcgs}.)  (In particular, then, $\Delta_0$ is the
usual spherical Laplacian on $S^2$.)  This was all noted by Newman and Penrose.\\

Of course $\{{}_sY_{lm}: -l \leq m \leq l\}$ is an orthonormal basis for ${\cal H}_{ls}$.  Note that
the map $f \rightarrow f^{R^{-1}} = R_s \circ f \circ R^{-1}$ is a {\em unitary mapping}
of $L^2({\bf L}^s)$ onto itself, for any $R \in SO(3)$.  Note also that for each such $R$,
${}_sY_{lm}^{R^{-1}}$ is an eigenfunction of $\Delta_s$ for the eigenvalue $\lambda_{ls}$.
(This follows from (\ref{edthgoR}) and (\ref{bedthgoR}).)  Thus, for any $R$,
$\{{}_sY_{lm}^{R^{-1}} : -l \leq m \leq l\}$ is also an orthonormal basis of ${\cal H}_{ls}$.
Note also that, by (\ref{rotsame}), we have, for $p \in U_R$,
\begin{equation}
\label{rotylm}
({}_sY_{lm}^{R^{-1}})_R(p) = {}_sY_{lmI}(R^{-1}p).
\end{equation}
and so, in the $(\theta_R,\phi_R)$ coordinates on $U_R$, we have
\begin{equation}
\label{sphaspnR}
({}_sY_{lm}^{R^{-1}})_R(\theta_R,\phi_R)
= e^{im\phi_R}(-1)^{m^{-}}\frac{a_{lm}}{[(l+s)!(l-s)!]^{1/2}}
\sin^{2l}(\theta_R/2)\sum_r \left( \displaystyle^{l-s}_{\:\:\:r}\right)\left( \displaystyle^{\:\:\:l+s}_{r+s-m}\right)
(-1)^{l-r-s}\cot^{2r+s-m}(\theta_R/2).
\end{equation}

\section{Zonal Harmonics and Kernels}

We now return to our new results.

Concerning the ${}_sY_{lm}$, we also make the following observation, which
is essentially in \cite{cc}, section 6.

\begin{proposition}
\label{yslmneval}
\[ \lim_{\theta \rightarrow 0^+} e^{is\phi} {}_sY_{lmI}(\theta,\phi) =
\begin{cases}
(-1)^{s^+} \left[\frac{2l+1}{4\pi}\right]^{1/2} & \mbox{if } m = -s\\
0 & \mbox{otherwise}
\end{cases}
\]
\end{proposition}
{\bf Proof} Since ${}_sY_{lmI}(\theta,\phi)$ equals $e^{im\phi}$ times a function of $\theta$,
this limit (which surely exists, by Theorem \ref{limns}) must be zero if $m \neq -s$.  If $m = -s$,
it follows from an examination of (\ref{sphareg}).  (In each term in the summation in (\ref{sphaspn}),
because of the binomial coefficients, we must have $r \leq l-s$.  As $\theta \rightarrow 0^+$,
$\sin^{2l}(\theta/2)\cot^{2r+s-m}(\theta/2) = \sin^{2l}(\theta/2)\cot^{2r+2s}(\theta/2)$
has a nonzero limit only if $2l \leq 2r+2s$.  These two conditions together allow us to look
only at the term in which $r=l-s$ when taking the limit, and the proposition now follows at once.)\\

We next discuss the {\em kernel} of the projection operator $P_{ls}$; we claim that it is a
{\em smooth-kernel operator}.  Since $P_{ls}$ maps sections in $L^2({\bf L}^s)$ to smooth sections
of ${\bf L}^s$, we must define what we mean by the kernel of such an operator.  For $x \in S^2$,
let ${\bf L}^s_x$ be $\pi^{-1}(x)$, the fiber above $x$.

\begin{definition}
\label{e1e2smtdf}
We say that $T: L^2({\bf L}^s) \rightarrow C^{\infty}({\bf L}^s)$ is a
{\em smooth-kernel operator} if $T$ has the form
\begin{equation}
\label{sopMdfwaybdl}
(Tf)(x) = \int_{S^2} K(x,y) f(y) dS(y)
\end{equation}
for any $f \in L^2({\bf L}^s)$; here we require:\\
(i) For each $(x,y) \in S^2 \times S^2$,
$K(x,y) : {\bf L}^s_y \rightarrow {\bf L}^s_x$ is linear;\\
(ii) $K(x,y)$ depends smoothly on $(x,y)
\in S^2 \times S^2$.

We also call the restriction of $T$ to $C^{\infty}({\bf L}^s)$
a smooth-kernel operator.
We call the collection of linear maps $K = \{K(x,y): (x,y)
\in S^2 \times S^2\}$ the {\em kernel} of $T$.
\end{definition}

Note here that the integral in (\ref{sopMdfwaybdl}) is a
vector-valued integral.  Also note that
the statement $(ii)$ means that the kernel is smooth in $(x,y)$
after one locally trivializes the bundle.

To say that $T$ is a smooth-kernel operator is evidently equivalent
to saying the following: for each $R, R' \in SO(3)$, there exists
a function $K_{R',R} \in C^{\infty}(U_{R'} \times U_R)$, such that
for all smooth sections $f$ of ${\bf L}^s$ with compact support in $U_R$,
we have
\begin{equation}
\label{krrprint}
(Tf)_{R'}(x) = \int_{S^2} K_{R',R}(x,y) f_R(y) dS(y),
\end{equation}
for all $x \in U_{R'}$.\\

Evidently, if $R, R', R_1, R_1' \in SO(3)$, if (\ref{krrprint}) holds, we must have

\begin{equation}
\label{krr1rprint}
K_{{R_1}',{R_1}}(x,y) = e^{is(\psi_1-\psi_2)}K_{R',R}(x,y)
\end{equation}

for all $x \in U_{{R_1}'} \cap U_{R'}$ and $y \in U_{{R_1}} \cap U_{R}$, where\\
\begin{equation}
\label{psi1}
\psi_1 \mbox{ is the angle from } \rho_{R'}(x) \mbox{ to } \rho_{{R_1}'}(x),
\mbox{ and } \psi_2 \mbox{ is the angle from } \rho_R(y) \mbox{ to } \rho_{R_1}(y).
\end{equation}

Let us now examine the kernel of the projection operator $P_{ls}$.  Evidently, if $f \in L^2({\bf L}^s)$,
we have
\begin{equation}
\label{plsfrst}
P_{ls}f = \sum_{m=-l}^l <f,\:{}_sY_{lm}>{}_sY_{lm}.
\end{equation}
(Here one could use ${}_sY_{lm}^{R_0^{-1}}$ in place of ${}_sY_{lm}$, for any $R_0 \in SO(3)$.)\\
Thus $P_{ls}$ has kernel $K^{ls}$, where for any $R, R' \in SO(3)$,
\begin{equation}
\label{kRdf}
K^{ls}_{R',R}(x,y) = \sum_{m=-l}^l {}_sY_{lmR'}(x)\overline{{}_sY_{lmR}(y)}.
\end{equation}
(Again here, one could use ${}_sY_{lm}^{R_0^{-1}}$ in place of ${}_sY_{lm}$, for any $R_0 \in SO(3)$.)\\

\begin{definition}
\label{szondf}
For $|l| \geq s$, we define ${}_sZ_l$, the ``$s$-zonal harmonic'' in ${\cal H}_{ls}$, by
\begin{equation}
\label{szondfway}
{}_sZ_l = (-1)^{s^+} \left[\frac{2l+1}{4\pi}\right]^{1/2}{}_{s}Y_{l,-s}.
\end{equation}
\end{definition}

We have the following new theorem, which provides justification for this definition.
When $s=0$ this theorem reduces to the
well-known result that
\begin{equation}
\label{zlokxy}
{}_0Z_l(x) := Z_l(x) = K^{l0}(x,{\bf N}) \mbox{ is the zonal harmonic in } {\cal H}_{l0}.
\end{equation}
\begin{theorem}
\label{plsker}
(a) Say ${\bf N} \in U_{R}$, and that $x \in U_I \cap U_{R'}$.  Then
\begin{equation}
\label{spnzonway}
{}_{s}Z_{lI}(x) =
e^{is(\psi_1-\psi_2)}K^{ls}_{R',R}(x,{\bf N}),
\end{equation}
where
\begin{equation}
\label{psiN}
\psi_1 \mbox{ is the angle from } \rho_{R'}(x) \mbox{ to } \rho_{I}(x),
\mbox{ and } \psi_2 \mbox{ is the angle from } \rho_R({\bf N}) \mbox{ to } \partial/\partial y
\mbox{ (at } {\bf N}).
\end{equation}
(b) For any $f \in {\cal H}_{ls}$, we have
\begin{equation}
\label{szonway}
\lim_{\theta \rightarrow 0} e^{is\phi}f_I(\theta,\phi) =
<f,\:{}_sZ_l>.
\end{equation}
\end{theorem}
{\bf Proof} For (a), note that if $y \in U_I \cap U_{R}$, then $K^{ls}_{I,I}(x,y)$ equals
\[ \sum_m {}_sY_{lmI}(x)\overline{{}_sY_{lmI}(y)} =
e^{is(\psi_1-\psi_2)}K^{ls}_{R',R}(x,y) \]
where
$\psi_1$ is the angle from $\rho_{R'}(x)$ to $\rho_{I}(x)$,
and $\psi_2$ is the angle from $\rho_R(y)$ to $\rho_{I}(y)$.
In this equation, let us say $y$ has coordinates $(\theta,\phi)$ (in the
usual coordinates on $U_I$).  Let us fix $\phi$, multiply both sides of the equation by
$e^{-is\phi}$, and take the limit as $\theta \rightarrow 0^+$.  (One can, if one
wishes, choose to take $\phi = 0$; or one can let $\phi$ be arbitrary, and recall our
choice of orientation at $\bf N$.)  (a) now follows at once from Proposition
\ref{yslmneval}.

To prove (b), one may assume $f = {}_sY_{lm}$ for some $m$; but then the result is evident
from Proposition \ref{yslmneval}.  (Alternatively, one can start with
the equation
$f_I(x) = \int_{S^2} K^{ls}_{I,I}(x,y) f_I(y) dS(y)$ (which is valid for
$f \in {\cal H}_{ls}$, $x \in U_I$), multiply by $e^{is\phi}$, and take a
limit.)  This completes the proof.\\
\ \\
For another new result, also familiar when $s = 0$, we examine $K^{ls}_{R,R}(x,x)$,
for $x \in U_R$.  It is evident from (\ref{krr1rprint}), (\ref{psi1}) that this quantity
is independent of $R$.  In fact it is independent of $x$ as well:

\begin{theorem}
\label{krrxx}
If $x \in U_R$, then
\begin{equation}
\label{krrxxway}
K^{ls}_{R,R}(x,x) = \frac{2l+1}{4\pi}.
\end{equation}
\end{theorem}
{\bf Proof} It is enough to prove this for $R=I$, since we will then know it unless $x = {\bf N}$ or
${\bf S}$, and the result will follow in those cases too by continuity.  Let us then take $R=I$ and
free the letter $R$ for other uses.

Say now $x,y \in U_I$.  Choose $R \in SO(3)$ with $Ry =x$.  Then $x \in U_I \cap U_R$, so
\[ K^{ls}_{I,I}(x,x) = K^{ls}_{R,R}(x,x) = \sum_m |({}_sY_{lm}^{R^{-1}})_R(x)|^2
 = \sum_m |{}_sY_{lmI}(R^{-1}x)|^2 = K^{ls}_{I,I}(y,y).\]
(We used (\ref{rotylm}).)  Accordingly $K^{ls}_{I,I}(x,x)$ is independent of $x \in SO(3)$.  Thus
for any fixed $x \in SO(3)$,
\[ 4\pi K^{ls}_{I,I}(x,x) = \int_{S^2} K^{ls}_{I,I}(y,y) dS(y) = \sum_m \int_{S^2} |{}_sY_{lm}(y)|^2 dS(y)
= 2l+1, \]
as desired.\\
\ \\
In doing estimates involving the $K^{ls}$, it will be very useful to observe the theorem
which follows.  In this theorem, $\edth^{[s]}_R$ will denote the differential operator on $U_R$ which
satisfies $(\edth^s F)_R = \edth^{[s]}_R F_R$, for any ordinary smooth function $F$ on $U_R$; explicitly,
$\edth^{[s]}_R = \edth_{s-1,R} \circ \ldots \edth_{0,R}$.  Similarly,
$\bedth^{[s]}_R$ will denote the differential operator on $U_R$ which
satisfies $(\bedth^s F)_R = \bedth^{[s]}_R F_R$, for any ordinary smooth function $F$ on $U_R$.

\begin{theorem}
\label{plskerest}
Say $x \in U_{R'}$ and $y \in U_R$.  Then if $s > 0$,
\begin{equation}
\label{klsested}
K^{ls}_{R',R}(x,y) =
\frac{(l-s)!}{(l+s)!}\edth^{[s]}_{R',x} \bedth^{[s]}_{R,y} K^{l0}(x,y),
\end{equation}
while if $s < 0$,
\begin{equation}
\label{klsestbed}
K^{ls}_{R',R}(x,y) =
\frac{(l+s)!}{(l-s)!}\bedth^{[-s]}_{R',x} \edth^{[-s]}_{R,y} K^{l0}(x,y).
\end{equation}
(Here $\edth^{[s]}_{R',x}$ denotes $\edth^{[s]}_{R'}$ taken in the $x$ variable, etc.)
\end{theorem}
{\bf Proof} This follows at once from (\ref{ylmblsed}), (\ref{ylmblsbed}),
and (\ref{kRdf}).\\

Of course $K^{l0}(x,y) = Z_{l,x}(y)$, the zonal harmonic in ${\cal H}_{l0}$
based at $x$.  It equals $[(2l+1)/4\pi]P_l^{1/2}(x\cdot y)$, where
$P^{\lambda}_l$ is the ultraspherical (or Gegenbauer)
polynomial of degree $l$ associated with $1/2$.

Series defining wavelets on the sphere, involving the $K^{l0}$, were estimated
in \cite{narc1}, \cite{narc2} and \cite{gmcw}.  Using these estimates, and Theorem
\ref{plskerest}, we will be able to estimate similar series, which involve the
$K^{ls}$ for general $s$, in place of the $K^{l0}$, in the next section.

\section{Spin Wavelets}

In sections 3-5 we have often used the variable ``$f$'' to denote a spin function on the
sphere.  In the remainder of this article, the variable ``$f$'' will be reserved for another purpose.

Specifically, say $f \in {\cal S}(\RR^+)$, $f \neq 0$,
and $f(0) = 0$.  Let  $({\bf M},g)$ be a smooth compact
oriented Riemannian manifold, and let $\Delta$ be the Laplace-Beltrami operator on
$\bf M$ (for instance, the spherical Laplacian if $\bf M$ is $S^2$).  Let $K_t$
be the kernel of $f(t^2\Delta)$.  Then, as we have explained in section 2, the functions
\begin{equation}
\label{manwavdf}
w_{t,x}(y) = \overline{K}_t(x,y),
\end{equation}
if multiplied by appropriate weights,
can be used as wavelets on ${\bf M}$.  In case ${\bf M} = S^n$ and $f$ has compact support
away from the origin, we shall say these wavelets are ``needlet-type''.

On the sphere $S^2$, we have
\begin{equation}
\label{sphkt}
K_t(x,y) = \sum_{l=0}^\infty \sum_{m=-l}^l f(t^2\lambda_l)  Y_{lm}(x) \overline{ Y_{lm}}(y).
\end{equation}

As we have explained in section 2, a key property of this kind of
wavelet is its {\em localization}.
Working in the general situation (general ${\bf M}$, general $f$), one has:\\
 for every pair of
$C^\infty$ differential operators $X$ (in $x$) and $Y$  (in $y$) on ${\bf M}$, and for every nonnegative integer
$N$, there
 exists $c:=c_{N,X,Y}$ such that for all $t>0$ and $x, y\in {\bf M}$
\begin{equation}
\label{diagloc}
\left| XYK_t(x,y) \right|\leq c ~ \frac{ t^{-(n+I+J)}}{(d(x,y)/t)^N},
\end{equation}
 where $I:=\deg X$ and $J:=\deg Y$.  In the case where ${\bf M}$ is the sphere, $f$ has compact support away from $0$,
and if one replaces $\lambda_l=l(l+1)$ by $l^2$ in the formula for $K_t$, this was
shown by Narcowich, Petrushev and Ward, in \cite{narc1} and \cite{narc2}.  The general result
was shown in \cite{gmcw}.  It will be important to note that in (\ref{diagloc}), the constant
$c$ can be chosen to be $\leq C\|f\|_M$, for some $M$.  (Here we have chosen a nondecreasing family of
norms $\{\|\:\|_M\}$ to define the Fr\'echet space topology on
${\cal S}_0(\RR^+) = \{f \in {\cal S}(\RR^+): f(0) = 0\}$.)  Indeed, this fact may be read off from the proof of
(\ref{diagloc}).  More simply, one may use the closed graph theorem: The collection of functions $K_t(x,y)$
(smooth in $t,x,y$ for $t > 0$) which satisfy (\ref{diagloc}) is itself naturally a Fr\'echet space, say $F$.
We want to see that the map taking $f$ to the function $K_{t,f}(x,y) \in F$ is continuous.
(Here we are temporarily letting $K_{t,f}$ be the kernel of $f(t^2\Delta)$.)   For this, by
the closed graph theorem, one need only note that the map from ${\cal S}_0(\RR^+)$ to $\CC$
which takes $f$ to $K_{t,f}(x,y)$ (for any fixed $t,x,y$) is continuous; but this is evident from
(\ref{sphkt}).\\

The proof of (\ref{diagloc}) in \cite{gmcw} would in fact go through without change if, instead of
$K_t$ being the kernel of $f(t^2\Delta)$, where $\Delta$ is the Laplace-Beltrami operator, one assumed
that $K_t$ was the kernel of $f(t^2D)$, where $D$ is any smooth positive second-order elliptic differential operator,
acting on $C^{\infty}({\bf M})$.  It is natural to conjecture then that something analogous
holds if $D$ is instead a smooth positive second-order elliptic differential operator, between sections
of line bundles on ${\bf M}$.

In fact, the following theorem says that this is the case for $D = \Delta_s: C^{\infty}({\bf L}^s) \rightarrow
C^{\infty}({\bf L}^s)$:

\begin{theorem}
\label{locestls}
Let $s$ be an integer.
Suppose $f \in {\cal S}(\RR^+)$, $f(0) = 0$.  Let $K_t$ be the kernel of $f(t^2\Delta_s)$.  Then:

For every $R, R' \in SO(3)$, every pair of compact sets ${\cal F}_R \subseteq U_R$ and ${\cal F}_{R'} \subseteq U_{R'}$,
and every pair of $C^{\infty}$ differential operators $X$ (in $x$) on $U_{R'}$ and $Y$ (in $y$) on $U_{R}$,
and for every nonnegative integer $N$, there exists $c$ such that for all $t > 0$, all $x \in {\cal F}_{R'}$
and all $y \in {\cal F}_{R}$, we have
\begin{equation}
\label{diaglocls}
\left| XYK_{t,R',R}(x,y) \right|\leq c ~ \frac{ t^{-(n+I+J)}}{(d(x,y)/t)^N},
\end{equation}
where $I:=\deg X$ and $J:=\deg Y$.
\end{theorem}
{\bf Note} The ``needlet-type'' case, where $f$ has compact support away from the origin, is easiest,
and that is the only case we will deal with in this article.  The general case will be dealt with
elsewhere. \\
{\bf Proof}
 We have
\begin{equation}
\label{flsexp}
f(t^2\Delta_s) = \sum_{l \geq |s|} f(t^2\lambda_{ls})P_{ls}
\end{equation}
(strong convergence), so that
\begin{equation}
\label{krexp}
K_{t,R',R}(x,y) = \sum_{l \geq |s|} f(t^2\lambda_{ls})K^{ls}_{R',R}.
\end{equation}
(The interchange of order of integration and summation, needed to derive (\ref{krexp}), is easily
justified through the rapid decay of $f$, and through (\ref{ylmpolygro}), once one recalls
(\ref{kRdf}).)
Now, say $s > 0$.  Then by Theorem \ref{plskerest},
\begin{equation}
\label{krexpgd}
K_{t,R',R}(x,y) = \edth^{[s]}_{R'x} \bedth^{[s]}_{Ry}\sum_{l \geq |s|} f(t^2\lambda_{ls})\frac{(l-s)!}{(l+s)!} K^{l0}(x,y),
\end{equation}
(The interchange of order of differentiation and summation, needed to derive (\ref{krexpgd}), is easily
justified through the rapid decay of $f$, and through (\ref{ylmpolygro}) in the case $s=0$, once
one recalls (\ref{kRdf}) in the case $s=0$.)  Now notice that for $n \in \ZZ$, we have
\begin{equation}
\label{lnnpl}
(l+n)(l+1-n) = l(l+1)-\gamma_n
\end{equation}
where
\[ \gamma_n = n(n-1). \]
Also note that
\[ \frac{(l+s)!}{(l-s)!} = [(l+s)(l+1-s)][[(l+s-1)(l+1-(s-1))]\ldots[(l+1)(l+1-1)]. \]
Accordingly
\begin{equation}
\label{lslsfac}
\frac{(l+s)!}{(l-s)!} = [l(l+1)-\gamma_s][l(l+1)-\gamma_{s-1}]\ldots[l(l+1)-\gamma_1].
\end{equation}
Note also that
\begin{equation}
\label{lamls}
\lambda_{ls} = l(l+1)-\gamma_{-s} = l(l+1)-s(s+1).
\end{equation}

To prove the theorem, let us first note that if we restrict to $t \in [T_0,T_1]$, where
$0 < T_0 < T_1$ are fixed, then the result is trivial (even if we do not assume $f$ has compact
support away from the origin).
Indeed, since $d(x,y) \leq 2\pi$ for all $x,y$, it is enough to show that $XYK_{t,R',R}$ is bounded
for $(t,x,y) \in [T_0,T_1] \times {\cal F}_{R'} \times {\cal F}_R$.  This, however, is evident from
(\ref{ylmpolygro}) and the rapid decay of $f$.\\

For other $t$, we shall assume that the support of $f$ is contained in some
interval $[a,b]$ with $b > a > 0$.  In that case, there is a $T_1 > 0$ such that
$K_{t,R,R'} \equiv 0$ for all $t > T_1$ (indeed, this will be true if $\lambda_{ss}T_1^2 > b$.)
Thus we need only prove that there exists a $T_0 > 0$ such that the result holds for all $t \in (0,T_0]$.
In fact, choose any $T_0 > 0$ with $\gamma_s T_0^2 < a/2$.
For $0 < t \leq T_0$, define
\begin{equation}
\label{ftdf}
f_t(u) = \frac{f(u-s(s+1)t^2)}{[u-\gamma_st^2][u-\gamma_{s-1}t^2]\ldots[u-\gamma_1t^2]}.
\end{equation}
The function $f_t$ is supported in $[a,b+s(s+1)T_0^2]$, a fixed compact interval.  (Note that the
denominator in (\ref{ftdf}) does not vanish for $u$ in this interval; in fact, each factor is
at least $a/2$.)   Then, in (\ref{krexpgd}), by (\ref{lslsfac}) we may write
\begin{equation}
\label{ftinest}
\sum_{l \geq s} f(t^2\lambda_{ls})\frac{(l-s)!}{(l+s)!} K^{l0}(x,y) =
t^{2s}\sum_{l \geq s} f_t(t^2\lambda_{l}) K^{l0}(x,y) = t^{2s}K_{[t]}(x,y)
\end{equation}
where $K_{[t]}$ is the kernel of $f_t(t^2\Delta)$.  Now, it is easy to see that
the functions $f_t$ ($t \in (0,T_0]$) form a
bounded subset of ${\cal S}_0(\RR^+)$.  Thus, by the remarks following (\ref{diagloc}),
 for every pair of
$C^\infty$ differential operators $X'$ (in $x$) and $Y'$  (in $y$) on $S^2$, and for every nonnegative integer
$N$, there exists $c$ such that whenever $0 < t \leq T_0$, and $x, y\in S^2$,
\begin{equation}
\label{diagloct}
\left| X'Y'K_{[t]}(x,y) \right|\leq c ~ \frac{ t^{-(n+I'+J')}}{(d(x,y)/t)^N},
\end{equation}
 where $I':=\deg X'$ and $J':=\deg Y'$.  In fact, this remains true if we insist $x \in {\cal F}_{R'}$,
$y \in {\cal F}_R$, and allow $X'$ (resp.\ $Y'$) to be smooth
differential operators only in a neighborhood of ${\cal F}_{R'}$
(resp.\ ${\cal F}_{R}$), since there are smooth differential
operators on all of $S^2$ which agree with these on ${\cal F}_{R'}$
(resp.\ ${\cal F}_{R}$). In the situation of (\ref{diaglocls}), note
that $\edth^{[s]}_{R'x}$ is a smooth differential operator of degree
$s$ in $x$ in a neighborhood of ${\cal F}_{R'}$, and
$\bedth^{[s]}_{Ry}$ is a smooth differential operator of degree $s$
in $y$ in a neighborhood of ${\cal F}_{R'}$.  Thus we have
\[\left| XYK_{t,R',R}(x,y) \right|
= t^{2s}\left|XY \edth^{[s]}_{R'x} \bedth^{[s]}_{Ry} K_{[t]}(x,y)\right|
\leq ct^{2s} ~ \frac{ t^{-[n+(I+s)+(J+s)]}}{(d(x,y)/t)^N}
=  c  \frac{ t^{-[n+I+J]}}{(d(x,y)/t)^N},\]
as desired.  Similar arguments work if $s < 0$.  This completes the proof.

\begin{remark}
\label{constsn}
\rm {As in the case $s=0$, the closed graph theorem implies that the constant $c$ in (\ref{diaglocls})
can be chosen to be $\leq C\|f\|_M$, for some $M$.  (Here, for $f \in {\cal S}(\RR^+)$, let
us define $\|f\|_M = \sum_{i+j \leq M} \|x^i \partial^j f\|$, where $\|\:\|$ denotes sup norm.)
 Also, $C$ and $M$ depend only on $I, J, N, R'$
and $R$.)  Note, then, that (\ref{diaglocls}) remains true if $K_t$ is the kernel of $f(t^2\Delta_s)$,
where now only $f \in {\cal S}_0^M(\RR^+) =
\{f \in C^M(\RR^+): f(0) = 0 \mbox{ and } \|f\|_M < \infty\}$.  This is evident, since ${\cal S}_0(\RR^+)
= \cap_L {\cal S}_0^L(\RR^+)$ is dense in ${\cal S}_0^M(\RR^+)$.  Again we will can choose $c$ to be
$\leq C\|f\|_M$, where $C$ and $M$ depend only on $I, J, N, R'$ and $R$.}
\end{remark}

Let us now explain what we shall mean by spin wavelets.  Say $f \in {\cal S}(\RR^+)$,
$f \neq 0$, $f(0)=0$.  Say $s$ is an integer, $R \in SO(3)$, $t > 0$, and $x \in U_R$.  Let us define
\begin{equation}
\label{needdf}
w_{txR} = \sum_{l \geq |s|} \sum_m \overline{f}(t^2\lambda_{ls})\:\overline{{}_sY_{lmR}}(x)\:{}_sY_{lm}.
\end{equation}
Then $w_{txR} \in C^{\infty}({\bf L}^s)$.  Also note that

\begin{equation}
\label{manwavspn}
w_{txR,R'}(y) = \overline{K}_{t,R,R'}(x,y)
\end{equation}
(if $x \in U_R$, $y \in U_{R'}$, and $K_t$ is the kernel of
$f(t^2\Delta_s)$).  This generalizes the case $s=0$ of (\ref{manwavdf}).  Moreover, if
$F \in L^2({\bf L}^s)$, then for $x \in U_R$
\begin{equation}
\label{betFdf}
\beta_{t,F,x,R} := \langle F, w_{txR}\rangle = (f(t^2\Delta_s)F)_R(x) = (\beta_{t,F})_R(x)
\end{equation}
where we have set $\beta_{t,F} = (f(t^2\Delta_s)F)$.  We call $\beta_{t,F,x,R}$ a {\em spin wavelet
coefficient} of $F$.  If $F = \sum_{l \geq |s|}\sum_m a_{lm}\:{}_sY_{lm} \in L^2({\bf L}^s)$, then
\begin{equation}
\label{betexp}
\beta_{t,F} = \sum_{l \geq |s|} \sum_m f(t^2\lambda_{ls})a_{lm}\:{}_sY_{lm}(x).
\end{equation}
so that
\begin{equation}
\label{betRexp}
\beta_{t,F,x,R} = \sum_{l \geq |s|} \sum_m f(t^2\lambda_{ls})a_{lm}\:{}_sY_{lmR}(x).
\end{equation}
In cosmological applications, one assumes that $f$ is {\em real-valued}.  In that case, by
(\ref{betexp} and (\ref{cinfsxpe}) -- (\ref{almmxp}), we have that
\begin{equation}
\label{spndem}
\beta_{t,F_{\e}} = [\beta_{t}(F)]_{\e};\:\:\:\beta_{t}(F_{\m}) = [\beta_{t}(F)]_{\m}
\end{equation}
for all $F \in L^2({\bf L}^s)$.\\
\ \\
Let us now explain how, in analogy to the case $s=0$, one can obtain a nearly tight frame
from spin wavelets.  Let $P = P_{|s|,s}$
be the projection onto the ${\cal H}_{|s|,s}$, the null space of $\Delta_{ls}$ (in $C^{\infty}({\bf L}^s)$.
Then for $a > 1$ sufficiently close to $1$,
for an appropriate discrete set  $\{x_{j,k}\}_{(j,k)\in \ZZ\times [1,N_j]}$ on $\bf M$,
certain $R_{j,k}$ with $x_{j,k} \in U_{R_{j,k}}$,
and certain weights $\mu_{j,k}$, the collection of $\{\phi_{j,k}:= \mu_{j,k}w_{a^j, x_{j,k}, R_{j,k}}\}$
constitutes a nearly tight wavelet frame for $(I-P)L^2({\bf L}^s)$.
(Note that this property is independent of the choice of
$R_{j,k}$.)  In that case, we will have that for some $C$, if $F \in (I-P)L^2({\bf L}^s)$, then, in $L^2$,
$F \sim C\sum_{j,k}\langle F,\phi_{j,k} \rangle \phi_{j,k}$ (note that this sum is {\em independent}
of the choice of $R_{j,k}$).  We need to make this all precise.\\

In order to find suitable $x_{j,k}$, we need to subdivide $S^2$ into a fine grid.  To this end,
we select $c_0, \delta_0 > 0$ with the following properties:\\
\ \\
(*) Whenever $0 < \tau < \delta$, we can write $S^2$ as a finite disjoint union of measurable sets
${\cal E}_k$ such that the diameter of each ${\cal E}_k$ is less than or equal to $\tau$, and
such that the measure of each ${\cal E}_k$ is at least $c_0\tau^2$.\\
\ \\
It is easy to see that such $\tau, \delta$ exist.  (We could even require that each set ${\cal E}_k$
be either a spherical cap, or a spherical rectangle bounded by latitude and longitude lines.  The
existence of $c_0, \delta$ as in (*) can even be established on general smooth compact oriented
Riemannian manifolds; see the discussion before Theorem 2.4 of \cite{gmfr}.)

By using Theorem \ref{locestls}, one can prove the following result, which generalizes Theorem 2.4
of \cite{gmfr} (when the manifold there is the sphere).

\begin{theorem}
\label{framainsp}
$(a)$ Fix $a >1$, and say $c_0 , \delta_0$ are as in (*) above.  Suppose $f \in {\mathcal S}(\RR^+)$,
and $f(0) = 0$.  Let $K_t$ be the kernel of $f(t^2\Delta_s)$, and let
$w_{txR}$ be as in (\ref{needdf}) and (\ref{manwavspn}) above.

Then there exists a constant $C_0 > 0$ $($depending only on $f, a, c_0$ and $\delta_0$$)$
as follows:\\
Let ${\mathcal J} \subseteq \ZZ$ be either finite, or equal to all of $\ZZ$.
Say $0 < b < 1$.
For each $j \in {\mathcal J}$, write $S^2$ as a finite disjoint union of measurable sets
$\{E_{j,k}: 1 \leq k \leq N_j\}$,
where:
\begin{equation}
\label{diamleq}
\mbox{the diameter of each } E_{j,k} \mbox{ is less than or equal to } ba^j,
\end{equation}
and where:
\begin{equation}
\label{measgeq}
\mbox{for each } j \mbox{ with } ba^j < \delta_0,\: \mu(E_{j,k}) \geq c_0(ba^j)^2.
\end{equation}
For each $j,k$, select $x_{j,k} \in E_{j,k}$ and select any $R_{j,k} \in SO(3)$ with $x_{j,k} \in
U_{R_{j,k}}$.   Set
\begin{equation}
\label{varpjkdf}
\varphi_{j,k} = w_{a^j, x_{j,k}, R_{j,k}} \in C^{\infty}({\bf L}^s).
\end{equation}
For $F \in L^2({\bf L}^s)$, set
\begin{equation}\notag
S^{\mathcal J}F :=
\sum_{j \in {\mathcal J}}\sum_k \mu(E_{j,k}) \langle  F, \varphi_{j,k}  \rangle  \varphi_{j,k}.
\end{equation}
Here, if ${\mathcal J} = \ZZ$, the series converges unconditionally in $L^2({\bf L}^s)$.

Let $Q^{\mathcal J} = \sum_{j \in {\mathcal J}} |f|^2(a^{2j} \Delta_s)$;
if ${\cal J} = \ZZ$, this series converges strongly on $L^2({\bf L}^s)$.
Then for all $F \in L^2({\bf L}^s)$,

\begin{equation}
\label{qsclose}
\left| \langle  (Q^{\mathcal J}-S^{\mathcal J})F,F  \rangle \right| \leq C_0b  \langle  F,F  \rangle
\end{equation}
$($or, equivalently, since $Q^{\mathcal J}-S^{\mathcal J}$ is self-adjoint,
$\|Q^{\mathcal J}-S^{\mathcal J}\| \leq C_0b$.$)$\\
$(b)$ In $(a)$, take ${\mathcal J} = \ZZ$; set $Q =
Q^{\ZZ}$, $S = S^{\ZZ}$.   Suppose that the Daubechies condition $(\ref{daub})$ holds.
Let $P = P_{|s|,s}$
be the projection onto the ${\cal H}_{|s|,s}$, the null space of $\Delta_{s}$ (in $C^{\infty}({\bf L}^s)$).
Then
\begin{equation}
\label{snrtgt}
(A_a-C_0b)(I-P) \leq S \leq (B_a + C_0b)(I-P)
\end{equation}
as operators on $L^2({\bf L}^s)$.  Thus, for any $F \in (I-P)L^2({\bf M})$,
\begin{equation}\notag
(A_a-C_0b)\|F\|^2 \leq \sum_{j,k} \mu(E_{j,k})| \langle  F,\varphi_{j,k}  \rangle |^2
\leq (B_a + C_0b)\|F\|^2,
\end{equation}
so that, if $A_a - C_0b > 0$, then
$\left\{ \mu(E_{j,k})^{1/2}\varphi_{j,k}\right\}_{j,k}$ is a frame for $(I-P)L^2({\bf L}^s)$,
with frame bounds $A_a - C_0b$ and $B_a + C_0b$.
\end{theorem}
{\bf Remarks}
1. By (\ref{daubest}) and (\ref{daubest2}), $B_a/A_a = 1 + O(|(a-1)^2 (\log|a-1|)|)$; evidently
$(B_a + C_0b)/(A_a - C_0b)$ can be made arbitrarily close to $B_a/A_a$ by
choosing $b$ sufficiently small.  Thus, Theorem \ref{framainsp} gives ``nearly tight" frames
for $(I-P)L^2({\bf L}^s)$.\\
2. The proof of Theorem \ref{framainsp} follows along the same lines of the proof for $s=0$, which is Theorem
2.4 of \cite{gmfr}.  It will be given elsewhere \cite{gms}.  In this article, we will only use Theorem
\ref{framainsp} in Section 8.\\
3. In the statement of Theorem 2.4 (a) of \cite{gmfr}, ${\cal J}$ was not assumed to be a finite set,
but rather, a cofinite set (that is, ${\cal J}^c$ was assumed to be finite), and one had unconditional convergence
of the sum for $S^{\cal J}F$ and strong convergence of the sum for $Q^{\cal J}$.  Note, however, that
this is equivalent to (a) for finite sets, since $Q^{\cal J} = Q^{\ZZ} - Q^{{\cal J}^c}$,
$S^{\cal J}F = S^{\ZZ}F - S^{{\cal J}^c}F$.

\section{Asymptotic Uncorrelation}

We define a random spin $s$ field $G(.)$ by assuming that there
exists a probability space $(\Omega ,\Im ,P)$ such that the map
$(x,\omega)\rightarrow G(x,\omega )$ is $\mathcal{B}(S^{2})\otimes
\Im $ measurable, $\Im $ denoting a $\sigma $-algebra on $\Omega $
and $\mathcal{B}(S^{2})$ the Borel $\sigma $-field of $S^{2}.$ 
(To clarify: $G(x,\omega)$ is a.e.\ in ${\bf L}^s_x$ (the fiber of ${\bf L}^s$ above $x$),
and to say that the map is measurable is to say that, for any $R \in SO(3)$
the complex-valued function over $U_R \times \Omega$ obtained from the map by trivializing the bundle 
over $U_R$ is measurable.)  As is customary, the dependence of $G$ upon 
$\omega$ will be suppressed in the notation.
Also, we assume that
\begin{equation*}
E\left[ \int_{S^{2}}\left| G(x)\right| ^{2}dS(x)\right] =C<\infty
\text{ ,}
\end{equation*}%
which in particular entails that $x\rightarrow G(x)$ belongs
to $L^{2}(\mathbf{L}^{s})$ with probability one (see for instance
\cite{partha}).

In cosmological applications it is natural to introduce an {\em
isotropy condition}:

\begin{definition}
\label{isot}
Let $G$ be a random spin $s$ field.  We say that $G$ is {\em isotropic} (in law)
if for every $x_1,\ldots,x_N \in U_I$, the joint probability distribution of
$G^R_I(x_1),\ldots,G^R_I(x_N)$ is independent of $R \in SO(3)$. \\
Similarly, if $G, H$ are both random spin $s$ fields defined on the
same probability space, we say that $G,H$ are {\em jointly
isotropic} (in law) if for every $x_1,\ldots,x_N,y_1,\ldots,y_M \in
U_I$, the joint probability distribution of
$G^R_I(x_1),\ldots,G^R_I(x_N), H^R_I(y_1),\ldots,H^R_I(y_M)$ is
independent of $R \in SO(3)$.
\end{definition}
Note that $G$ is isotropic if and only if, for every $R' \in SO(3)$, and for every
$p_1,\ldots,p_N \in U_{R'}$, the joint probability distribution of
$G^R_{R'}(p_1),\ldots,G^R_{R'}(p_N)$ is independent of $R \in SO(3)$.
Indeed, say $p_j = R'x_j$,
$x_j \in U_I$ (for $j=1,\ldots,N$).   One need only then observe that, by (\ref{rotsame}),
$G^R_{R'}(p_j) =  G^{RR'}_I(x_j)$ for each $j$.  Similarly, $G,H$ are jointly isotropic
if and only if, for every $R' \in SO(3)$, and for every
$p_1,\ldots,p_N,q_1,\ldots,q_M \in U_{R'}$, the joint probability distribution of
$G^R_{R'}(p_1),\ldots,G^R_{R'}(p_N),G^R_{R'}(q_1),\ldots,G^R_{R'}(q_N)$ is independent of $R \in SO(3)$.\\

If $G$ is isotropic, then for any $R' \in SO(3)$ and
any $p \in U_{R'}$, the probability distribution of $G^R_{R'}(p)$ is independent of
$R \in SO(3)$.  Restricting $R$ to the subgroup of $SO(3)$ which fixes
$p$, we see in particular that $G_{R'}(p)$ has the same probability distribution as
$G_{R'R}(p) = e^{is\psi} G_{R'}(p)$, where $\psi$ is the angle from $\rho_{R'}(p)$ to $\rho_{R'R}(p)$.
If $s \neq 0$, we may take $\psi$ with $e^{is\psi} = -1$.  We then see that the probability distributions
of $\Re G_{R'}(p)$ is the same
as that of $-\Re G_{R'}(p)$, similarly for $\Im G_{R'}(p)$; in particular, both
$\Re G_{R'}(p)$ and $\Im G_{R'}(p)$ must have expectation $0$.    Again if $s \neq 0$,
we may take $\psi$ with $e^{is\psi} = i$, and we
see now that $\Re G_{R'}(p)$ and $\Im G_{R'}(p)$ have the same probability distributions.\\

Say $G$ is isotropic, and that $G^R$ has spin $s$ spherical harmonic
coefficients $a_{lm}^R$ (for $R \in SO(3)$).  Then, evidently, the
$a_{lm}^R$ all have the same probability distribution, and each $a_{lm}^R$ has expectation zero.
(Here $a_{lm}^R = \langle G^R,\:_sY_{lm} \rangle
= \langle G^R_{R'},\:_sY_{lmR'} \rangle$ (for any $R' \in SO(3)$) where $\langle \: \rangle$ is the
inner product on $S^2$.).  \\

If we take $R$ to be a rotation about ${\bf N}$, we see (from
$a_{lm}^R = \langle G,\:_sY_{lm}^{R^{-1}} \rangle$) and from (\ref{sphaspn}) that, if $m \neq 0$,
then the probability distribution of $a_{lm}$ is the same as that of $e^{i\theta}a_{lm}$ for any
$\theta$.  In particular, $\Re a_{lm}$ and $\Im a_{lm}$ will have the same probability distributions,
if $m \neq 0$.\\
If $G, H$ are jointly isotropic, $a_{lm}^R = \langle G^R,\:_sY_{lm} \rangle$,
$b_{lm}^R = \langle H^R,\:_sY_{lm} \rangle$, then for any
$l,m,l',m'$ the covariances $E(a_{lm}^R\overline{b}_{l'm'}^R)$ are evidently
independent of $R$.\\

In fact, as in the case $s=0$ of \cite{bmv07}, we can say much more
about these covariances:
\begin{theorem}
\label{schurthm}
Suppose $G,H$ are jointly isotropic random spin $s$ fields.  Let $a_{lm} = \langle G^R,\:_sY_{lm} \rangle$,
and $b_{lm} = \langle H^R,\:_sY_{lm} \rangle$.  Then
$E(a_{lm}\overline{b}_{l'm'}) = 0$ unless $l=l'$ and $m=m'$.  Moreover
$E(a_{lm}\overline{b}_{lm})$ does not depend on $m$; we will denote it by $C_{l,G,H}$.
\end{theorem}
{\bf Proof}
We define {\em Wigner's $D^{l}$ matrices} by
\begin{equation}
 _{s}Y_{lm}^R =\sum_{m}D_{m^{\prime }m}^{l}(R)
\:_{s}Y_{lm^{\prime }} \text{ .}  \label{wig}
\end{equation}
These matrices surely exist and are unitary, since the $_sY_{lm}^R$
form an orthonormal basis of ${\cal H}_{ls}$ whenever $R \in SO(3)$
and $l \geq |s|$.  The matrices are {\em independent} of $s$, as one
sees at once from (\ref{wig}), use of the $\edth$ and $\bedth$
operators, and (\ref{ylmblsed}), (\ref{ylmblsbed}), (\ref{edthgoR}),
(\ref{bedthgoR}).  Since $ _{s}Y_{lm}^{RR'} =\sum_{m}D_{m^{\prime
}m}^{l}(R) \:_{s}Y_{lm^{\prime }}^{R'}$ for any $R, R'$, we see that
in fact the map $R \to D_{m^{\prime }m}^{l}(R)$ is a unitary
representation of $SO(3)$, see also \cite{vmk} and \cite{vilkli}.
This representation is evidently isomorphic to the action of $SO(3)$
on ${\cal H}_{l0}$, which as is well known, is irreducible;
moreover, the representations for different $l$ are not equivalent.

Now let $a_{lm}^R = \langle G^R,\:_sY_{lm} \rangle$, and $b_{lm}^R = \langle H^R,\:_sY_{lm} \rangle$.
For each $l, l'$ define the covariance matrix $M^{ll'}$ by $M^{ll'}_{mm'} = E(a_{lm}\overline{b}_{l'm'})
= E(a_{lm}^R \overline{b}_{l'm'}^R)$ for any $R$.  Noting that
$a_{lm}^R = \langle G,\:_sY_{lm}^{R^{-1}} \rangle$, one easily computes that
$M^{ll'} = [D^l(R)]^*M^{ll'}D^l(R)$ for any $R$.  The theorem now follows at once from Schur's lemma.\\

In cosmological applications, having an``asymptotic uncorrelation'' theorem
of the following kind is of great importance.  Generalizing a key result in \cite{BKMP06-3}, we now show:

\begin{theorem}
\label{spinuncor}
Let $G,H$ be random, jointly isotropic spin $s$ fields.  If $x \in U_R$, let $w_{txR}$ be the
spin wavelet of (\ref{needdf}),
where $f$ has compact support away from the origin.
Assume that $C_{l,G,H} = g(l)$ for a smooth function $g$ on
the interval $(|s|,\infty)$, which satisfies
the following condition, for some $\alpha > 2$:
for every $i \in \NN_0$ there exists $c_i > 0$ such that
\begin{equation}
\label{gelest}
|g^{(i)}(u)| \leq c_iu^{-\alpha-i}.
\end{equation}
Assume also that, for some $c > 0$, certain $\rho, \tau \in \RR$, and for all sufficiently large $l$,
\begin{equation}
\label{clglow}
C_{l,G,G} \geq c l^{-\rho}, \mbox{   and   }  C_{l,H,H} \geq c l^{-\tau}.
\end{equation}
For $x \in U_R$, let $\beta_{G,t,x,R} = <G,w_{txR}>$, $\beta_{H,t,x,R} = <H,w_{txR}>$.
Then for any $M \in \NN_0$, there exists $T_0 > 0$ and a constant ${\cal C}_M > 0$, such that for all $x,y \in U_R$,
$0 < t < T_0$, we have
\begin{equation}
\label{spinuncorway}
|\mbox{Cor}(\beta_{G,t,x,R},\beta_{H,t,y,R})| \leq t^{\alpha-(\rho+\tau)/2}\frac{{\cal C}_M}{(d(x,y)/t)^M}.
\end{equation}
In particular, for fixed $x,y$, $|\mbox{Cor}(\beta_{G,t,x,R},\beta_{H,t,y,R})| \to 0$ as $t \to 0^+$.
\end{theorem}
{\bf Remark}: Note that $|\mbox{Cor}(\beta_{t,x,R},\beta_{t,y,R})|$ is independent of $R$.\\
{\bf Proof}
By definition,
\begin{equation}
\label{cordf}
|\mbox{Cor}(\beta_{G,t,x,R},\beta_{H,t,y,R})| =
\frac{|E(\beta_{G,t,x,R}\overline{\beta}_{H,t,y,R})|}{\sqrt{E(|\beta_{G,t,y,R}|^2)}\sqrt{E(|\beta_{H,t,y,R}|^2)}}.
\end{equation}
Say $G, H$ have spin $s$ spherical harmonic coefficients $a_{lm}, b_{lm}$
respectively.  Note that
\begin{equation}
\label{betexp1}
\beta_{G,t,x,R} = \sum_{l \geq |s|} \sum_m f(t^2\lambda_{ls})a_{lm}\:{}_sY_{lmR}(x),
\end{equation}
and similarly for $H$.  Accordingly
\begin{eqnarray}
E(\beta_{G,t,x,R}\overline{\beta}_{H,t,y,R}) & = & \sum_{l \geq |s|}
|f|^2(t^2\lambda_{ls})C_{l,G,H}\sum_m {}_sY_{lmR}(x)\overline{{}_sY_{lmR}}(y)  \\
& = & \sum_{l \geq |s|} |f|^2(t^2\lambda_{ls})C_{l,G,H}K^{ls}_{R,R}(x,y) \label{corkls}.
\end{eqnarray}
To estimate $E(|\beta_{G,t,x,R}|^2)$, let us note that, since $f \neq 0$, we may choose $c_0$ and
$0 < a < b$ such that $|f|^2 \geq c_0$ on $[a^2,b^2]$.  Thus
$|f|^2(t^2\lambda_{ls}) \geq c_0$ if $a/t \leq \sqrt{\lambda_{ls}} \leq b/t$.  Choose $a',b'$ with
$a < a' < b' < b$.  Since $\lim_{l \rightarrow \infty} l/\sqrt{\lambda_{ls}} = 1$, there is a $T_1 > 0$
such that if $0 < t < T_1$, $|f|^2(t^2\lambda_{ls}) \geq c_0$ if $a'/t \leq l \leq b'/t$.  Choose
$0 < T_2 < T_1$ with $(b'-a')/t > 1$ if $0 < t < T_2$.  Then, for $t$ sufficiently small, by Theorem \ref{krrxx}
and (\ref{clglow}), for some $C', C > 0$,
\begin{equation}
\label{extrest}
E(|\beta_{G,t,x,R}|^2) \geq \sum_{\frac{a'}{t} \leq l \leq \frac{b'}{t}} c_0 l^{-\rho}\frac{2l+1}{4\pi}
\geq C'\frac{b'-a'}{t} t^{\alpha} t^{-1} = Ct^{-2+\rho}.
\end{equation}
Similarly,
\begin{equation}
\label{extresth}
E(|\beta_{H,t,x,R}|^2) \geq Ct^{-2+\tau}.
\end{equation}

We estimate the numerator in (\ref{cordf}) by using Theorem \ref{locestls}.  In order
to do this, we must transform $g$ in some simple ways.

Say $s \geq 0$.  Recall that $\lambda_{ls} = (l-s)(l+s+1) = l(l+1)-s(s+1)$.  Thus, if $v > 0$, $v = \lambda_{ls}$ if
and only if $l =  \frac{1}{2}[-1+\sqrt{1+4[v+s(s+1)}] := Q(v)$, say.
Then $Q: (0,\infty) \rightarrow (0,\infty)$ bijectively.  For $v > 0$, set $h(v) = g(Q(v))$;
then $g(l) = h(\lambda_{ls})$.  Moreover, it is easy to see, from (\ref{gelest}), that for every $i \in \NN_0$
there exists $C_i' > 0$
such that
\begin{equation}
\label{hsymbest}
|\frac{d^i}{dv^i}h(v)| \leq C_i'v^{-\alpha/2-i}
\end{equation}

From (\ref{corkls}), we have
\begin{eqnarray}
E(\beta_{txR}\overline{\beta}_{tyR})
& = & \sum_{l \geq |s|} |f|^2(t^2\lambda_{ls})h(\lambda_{ls})K^{ls}_{R,R}(x,y) \\
& = & t^{\alpha}K_{t,R,R}(x,y),
\end{eqnarray}
where $K_t$ is the kernel of $f_{[t]}(t^2\Delta_s)$, where now
\begin{equation}
\label{ftvdf}
f_{[t]}(w) = |f|^2(w)[t^{-\alpha}h(w/t^2)],
\end{equation}
(since $|f|^2(t^2\lambda_{ls})h(\lambda_{ls}) = t^{\alpha}f_{[t]}(t^2\lambda_{ls})$).
Select $0 < A < B$ with supp$f \subseteq [A,B]$.  Then supp$f_{[t]} \subseteq [A,B]$
for all $t$.  By use of (\ref{hsymbest}), one sees readily that for any $L \in \NN_0$,
there is a $C > 0$ such that $\|f_{[t]}\|_{C^L} \leq C$ for all $t$.  Accordingly,
by Theorem \ref{locestls} and Remark \ref{constsn}, for every $M$ there exists ${\cal C}_M' > 0$ with
\[ |E(\beta_{txR}\overline{\beta}_{tyR})| \leq \frac{{\cal C}_M't^{\alpha}t^{-2}}{(d(x,y)/t)^M}\]
Using this in (\ref{cordf}) together with (\ref{extrest}), we find the desired result.
Similarly for $s < 0$.\\
\ \\
{\bf Remarks} 1. Since supp$f_{[t]} \subseteq [A,B]$, when deriving the estimate $\|f_{[t]}\|_{C^L} \leq C$
from (\ref{hsymbest}) for $t$ small, one need only use (\ref{hsymbest}) for $v$ large (by (\ref{ftvdf})).
Thus, instead of assuming $g$ is smooth
on $(|s|,\infty)$, and that (\ref{gelest}) holds there, one may assume only that for some $T > 0$,
$g$ is smooth on $(T,\infty)$, and that (\ref{gelest}) holds there.  Equivalently, we are assuming
that $g(u)$ agrees with an ordinary symbol of order $-\alpha$ for large $u$.\\
2. By remark (\ref{constsn}), an examination of the proof of the theorem shows that for each $M$
there is an $L \in \NN_0$ such that, instead of assuming that $g$ is smooth
and satisfies (\ref{gelest}) for all $i$, we need only assume that $g \in C^L$ and
satisfies (\ref{gelest}) for all $i \leq L$.\\
3. The hypothesis (\ref{gelest}), which is used only to estimate the derivatives of
$h$ and then the derivatives of $f_{[t]}$, can evidently be relaxed -- at a cost -- to
$|g^{(i)}(u)| \leq c_iu^{-\alpha_i}$ whenever $i \leq L$, for certain $\alpha_i$.  The cost in doing this
is that one would need to multiply the right side of (\ref{spinuncorway}) by a factor of $t^{-N}$ for some
sufficiently large $N$.  It would still evidently follow that
for fixed $x,y$, $|\mbox{Cor}(\beta_{G,t,x,R},\beta_{H,t,y,R})| \to 0$ as $t \to 0^+$.\\
4. If $G=H$, in applications of the theorem, it is natural to assume $\alpha = \rho = \tau$. \\
5. In cosmological applications, where we have (\ref{spndem}), and where $s=2$, one assumes that
the polarization field is (a single sample of) an isotropic random field $F$, and that
$F_{\e}$ and $F_{\m}$ are jointly isotropic.

\section{Stochastic Limit Theorems}

In this final section, we give a brief glimpse of the significance
of the results of the previous section in cosmology.   We prove some
stochastic limit results, and then indicate why they are important.

Say $\epsilon > 0$.

We use a nearly tight frame, as provided by Theorem \ref{framainsp}.  Specifically, in that theorem, we
fix a real $f$ supported in $[1/a^2,a^2]$ for which the Daubechies sum
$\sum_{j=-\infty}^{\infty} f^2(a^{2j} u) = 1$ for all $u > 0$, so that $A_a = B_a = 1$.
We then choose $b$ sufficiently small that $C_0b < \epsilon$.  We then produce $\{\phi_{j,k}\}$
as in Theorem \ref{framainsp}, and let $\phi_{j,k} = \mu(E_{j,k})^{1/2}\varphi_{j,k}$, so that
$\{\phi_{j,k}\}$ is a nearly tight frame for $(I-P)L^2({\bf L}^s)$, with frame bounds $1-\epsilon$
and $1 + \epsilon$.  Also, we have (\ref{qsclose}) for
general finite subsets ${\cal J}$.  Of course our choice of $\{\phi_{j,k}\}$ depends on $\epsilon$;
if we want to indicate the dependence on $\epsilon$, we will write $\phi_{j,k,\epsilon}$.

Say $G = \sum A_{lm}\:_sY_{lm}$
is an isotropic random spin $s$ field, and let $C_l = C_{l,G,G}$ (notation as in Theorem
\ref{spinuncor}).   For $j \in \ZZ$, we also let
\begin{equation}
\label{ejdf}
{\gamma}_j = \sum_{a^{-2} \leq a^{2j}\lambda_{ls} \leq a^{2}} C_l(2l+1).
\end{equation}
(To avoid confusion, the summation in (\ref{ejdf}), and similar summations in the sequel,
are over those $l$ for which $\lambda_{ls}$ is defined and
for which $a^{-2} \leq a^{2j}\lambda_{ls} \leq a^{2}$; recall that $\lambda_{ls}$ is defined for $l \geq |s|$.)
Note that, for any $l$, there at most $3$ values of $j$ for which $a^{-2} \leq a^{2j}\lambda_{ls} \leq a^{2}$,
so that
\begin{equation}
\label{gamvar}
\sum_j \gamma_j \leq 3\sum_{l \geq |s|} C_l(2l+1) = 3\mbox{Var}G.
\end{equation}

We let $\beta_{jk} = \beta_{jk\epsilon} = \langle G,\phi_{j,k,\epsilon}\rangle$.  We focus on the quadratic statistics

\begin{equation}
\label{gammatdf}
\tilde{\Gamma}_{j} = \sum_{k}\left| \beta_{jk}\right| ^{2}
\end{equation}
and
\begin{equation}
\label{gammadf}
\widehat{\Gamma}_{j} = \sum_{l \geq |s|,\:m}f^{2}(a^{2j}\lambda _{ls})|A_{lm}|^2.
\end{equation}
Evidently
\begin{equation}
\label{gammaxpdf}
E(\widehat{\Gamma}_{j}) = \sum_{l \geq |s|} f^{2}(a^{2j}\lambda _{ls})C_l(2l+1) =
\sum_{a^{-2} \leq a^{2j}\lambda_{ls} \leq a^{2}} f^{2}(a^{2j}\lambda _{ls})C_l(2l+1).
\end{equation}

We have:

\begin{proposition}
\label{egam}
$E(|\widehat{\Gamma}_{j}- \tilde{\Gamma}_{j}|) \leq \epsilon \gamma_j$
for all $j$.
\end{proposition}
{\bf Proof}  For each $j \in \ZZ$, let
\begin{equation}
\label{gjdef}
G_j = \sum_{a^{-2} \leq a^{2j}\lambda_{ls} \leq a^{2}} A_{lm}\:_sY_{lm}.
\end{equation}
Then
\[ \beta_{jk} = \langle G,\phi_{j,k} \rangle =
\mu(E_{j,k})^{1/2}\sum_{l \geq |s|} \sum_m f(a^{2j}\lambda_{ls})A_{lm}\:{}_sY_{lm} =
\langle G_j,\phi_{j,k} \rangle. \]
Note also that $(I-P)G_j = G_j$, since no term with $l = |s|$ can appear in the summation
in (\ref{gjdef}) (for then $\lambda_{ls}=0$).
By Theorem \ref{framainsp} (a) with ${\cal J} = \{j\}$, we see that, in the notation of
that theorem,
\[ \tilde{\Gamma}_j = \langle S^{\cal J}G_j,G_j \rangle, \]
\[ \widehat{\Gamma}_{j} = \langle Q^{\cal J}G_j,G_j \rangle, \]
so that, by Theorem \ref{framainsp} (a),
\[ |\widehat{\Gamma}_{j}- \tilde{\Gamma}_{j}| \leq \epsilon\|G_j\|_2^2 =
\epsilon \sum _{a^{-2} \leq a^{2j}\lambda_{ls} \leq
a^{2}}|A_{lm}|^2. \]
The proposition now follows at once if we take expectations of both sides.\\

In Proposition \ref{egam}, $\tilde{\Gamma}_{j} = \sum_{k}\left|
\beta_{jk\epsilon}\right| ^{2}$ depends on our choice of $\epsilon$,
and we can even let $\epsilon << 1$ depend on $j$ here.  By
Tchebychev's inequality, the probability distribution of
$\tilde{\Gamma}_{j}$ is then a small perturbation of
that of $\widehat{\Gamma}_{j}$.  We conclude by briefly discussing the latter. \\

We will look at random spin $s$ fields $G = \sum A_{lm}\:_sY_{lm}$ for which
$\overline{A}_{lm} = A_{l,-m}$ for all $m$ (so that, in particular, $A_{l0}$ is real).
For want of a better word, let us
call such a $G$ {\em involutive}.  Note that, by (\ref{cinfsxpe}) -- (\ref{almmxp}), if $F$
is a random spin $s$ field, then both $F_{\e}$ and $F_{\m}$ are involutive.

We will also be assuming that $G$ is {\em Gaussian}, by which we mean that
$\{\Re G_I(x): x \in U_I\} \cup \{\Im G_I(x): x \in U_I\}$ is a Gaussian family.
In that case ${\cal F} := \{\Re A_{lm}\} \cup \{\Im A_{lm}\}$ is also a Gaussian family.
Recall that, if $m \neq 0$, $\{\Re A_{lm}\}$ and $\{\Im A_{lm}\}$ have the same
probability distributions, so that they have variances $E\left| \Re A_{lm}\right|^2,\:
E\left| \Im A_{lm}\right| ^{2}=C_{l}/2$.  Since we are assuming $G$ is involutive,
$\Re A_{lm} = (A_{lm} + A_{l,-m})/2$, $\Im A_{lm} = -i(A_{lm} - A_{l,-m})/2$.  It follows easily
from this and
from Theorem \ref{schurthm} that the elements of ${\cal F}$ are pairwise uncorrelated.
Since ${\cal F}$ is a Gaussian family they are independent.

We then have:

\begin{proposition}
\label{centlim}
For all $j=1,2,...,$,
\begin{eqnarray*}
Var\left\{ \widehat{\Gamma }_{j}\right\} &=&E\left[ \widehat{\Gamma }%
_{j}-E\widehat{\Gamma }_{j}\right] ^{2}=\sum_{l \geq |s|}2f^{4}(a^{2j}\lambda
_{ls})(C_{l})^{2}(2l+1)\text{ .%
}
\end{eqnarray*}%
Say now that for some $c, C > 0$ and some $\alpha > 2$, $cl^{-\alpha} \leq C_l \leq Cl^{-\alpha}$
for $l > 0$.  Then
\begin{equation*}
\frac{\widehat{\Gamma }_{j}-E\widehat{\Gamma }_{j}}{\sqrt{Var\left\{
\widehat{\Gamma }_{j}\right\} }}\rightarrow _{d}N(0,1)\text{ as }%
j\rightarrow -\infty \text{ ,}
\end{equation*}%
$\rightarrow _{d}$ denoting as usual convergence in probability law, and
$N(0,1)$ denoting as usual the standard normal distribution.
\end{proposition}
\textbf{Proof }
Write%
\begin{equation*}
X_{lm}(j):=\left\{
\begin{array}{c}
\sqrt{2}f(a^{2j}\lambda_{ls})\Re A_{lm}
\text{ for }m>0 \\
f(a^{2j}\lambda_{ls})A_{l0} \text{ for }m =0 \\
\sqrt{2}f(a^{2j}\lambda_{ls})\Im A_{lm} \text{ for }m<0. \\
\end{array}%
\right. ;
\end{equation*}%
For each $j$, $\left\{ X_{lm}(j)\right\} $ is a triangular array of
independent heteroscedastic (\emph{= unequal variance}) Gaussian
random variables such that
\begin{eqnarray*}
\widehat{\Gamma }_{j} &=&\sum_{l \geq |s|,\: m}X_{lm}^{2}(j)\text{ , }%
EX_{lm}^{2}(j)=f^{2}(a^{2j}\lambda_{ls})C_{l}%
\text{ ,} \\
Var\left\{ X_{lm}^{2}(j)\right\} &=&EX_{lm}^{4}(j)-\left\{
EX_{lm}^{2}(j)\right\} ^{2}=2f^{4}(a^{2j}\lambda_{ls}%
)(C_{l})^{2}\text{ .}
\end{eqnarray*}%
Hence we have easily%
\begin{eqnarray*}
Var\left\{ \widehat{\Gamma }_{j}\right\} &=&\sum_{l \geq |s|,\:m}Var\left\{
X_{lm}^{2}(j)\right\} =\sum_{l \geq |s|,\:m}2f^{4}(a^{2j}\lambda_{ls}%
)(C_{l})^{2} \\
&=&\sum_{l \geq |s|}2f^{4}(a^{2j}\lambda_{ls})(C_{l})^{2}(2l+1)%
\text{ .}
\end{eqnarray*}%
To establish the Central Limit Theorem, it is then enough to check
the validity of the Lindeberg-Levy condition (\cite{davidson}) which
here takes the simple form%
\begin{eqnarray*}
&&\lim_{j\rightarrow -\infty }\frac{\max_{a^{-2} \leq a^{2j}\lambda_{ls} \leq a^{2}}%
\max_{m}Var(X_{lm}^{2})}{\sum_{a^{-2} \leq a^{2j}\lambda_{ls} \leq a^{2}}2f^{4}(a^{2j}\lambda_{ls}%
)(C_{l})^{2}} \\
&\leq &K\lim_{j\rightarrow -\infty }\frac{\max_{a^{-2} \leq a^{2j}\lambda_{ls} \leq a^{2}}
\max_{m}f^{4}(a^{2j}\lambda_{ls}%
)(C_{l})^{2}}{\sum_{a^{-2} \leq a^{2j}\lambda_{ls} \leq a^{2}}
f^{4}(a^{2j}\lambda_{ls}%
)(C_{l})^{2}(2l+1)} \\
&\leq &K\lim_{j\rightarrow -\infty }\frac{a^{2j\alpha }}{a^{2j\alpha}
\sum_{a^{-2} \leq a^{2j}\lambda_{ls} \leq a^{2}} f^{4}(a^{2j}\lambda_{ls})l}=0\ .
\end{eqnarray*}
This completes the proof.

\end{document}